\documentclass[a4paper,11pt]{article}
\usepackage[active]{srcltx}
\usepackage{amsmath}\usepackage{amssymb}\usepackage{amscd}

\usepackage[]{graphicx}

\newcommand{\tr}{\tilde{r}}
\newcommand{\tR}{\tilde{R}}
\newcommand{\eC}{1_{{}_{C_2}}}
\newcommand{\eCn}{1_{{}_{C_2^n}}}
\newcommand{\epsCn}{\epsilon_{{}_{C_2^n}}}
\newcommand{\eA}{1_{{}_{A_{n-1}}}}

\newcounter{pic}

\newenvironment{pic}[1][\bf Fig. \arabic{pic}]{
        
        \refstepcounter{pic}\noindent\textbf{#1.}${}$\hspace{5pt}${}$\it}{}

\newtheorem{defi}{Definition}[section]

\newtheorem{prop}[defi]{Proposition}

\usepackage{empheq}

\begin{document}

\begin{center}

\bigskip
{\Large\bf Alternating subgroups of Coxeter groups and their spinor extensions}

\vspace{1cm}

{\large {\bf O. V. Ogievetsky$^{\circ\diamond}$\footnote{On leave of absence from P. N. Lebedev Physical Institute, Leninsky Pr. 53,
117924 Moscow, Russia} and L. Poulain d'Andecy$^{\circ}$}}

\vskip 0.6cm

$\circ\ ${\large Center of Theoretical Physics, Luminy \\
13288 Marseille, France}

\vspace{.6cm}
$\diamond\ ${\large 
J.-V. Poncelet French-Russian Laboratory, UMI 2615 du CNRS, Independent University of Moscow, 11 B. Vlasievski per., 119002 Moscow, Russia}

\end{center}

\vskip 1.2cm
\begin{abstract}\noindent
Let $G$ be a discrete Coxeter group, $G^+$ its alternating subgroup
and $\tilde{G}^+$ the spinor cover of $G^+$. A presentation of the groups $G^+$ and $\tilde{G}^+$ is proved for an arbitrary Coxeter system $(G,S)$; the generators are related to edges of the Coxeter graph. Results of the Coxeter--Todd algorithm - with this presentation - for the chains of alternating groups of types A, B and D are given.

\end{abstract}

\section{{\hspace{-0.55cm}.\hspace{0.55cm}}Introduction}

Let $(G,S)$ be a Coxeter system and $G^+$ the alternating subgroup of $G$. A presentation of $G^+$ by generators and relations  is given, as an 
exercise, in \cite{Bour}. One (arbitrary) vertex in the Coxeter graph is distinguished in this presentation. 

\vskip .2cm
For the Coxeter groups of type $A$ (that is, symmetric groups) the standard Artin presentation equips the chain of symmetric groups with a structure of a local and stationary chain of groups. The representation theory of the chain of symmetric groups was revisited in \cite{OV}.
It was emphasized in \cite{OV} that the method developed there for the representation theory of the symmetric groups could be hopefully relevant for a wider class of local and stationary chains of groups or algebras   
(see \cite{Vershik91} for the definition of locality and stationarity). For type A, the Bourbaki presentation of the chain of 
the alternating groups is not local. 

\vskip .2cm
A discrete Coxeter group $G$ is a subgroup of an orthogonal group $O$; the spin cover of $O$ restricts to a central extension $\tilde{G}$ of $G$.
A presentation of $\tilde{G}$ can be found in \cite{Mo2}.
The subgroup $G^+$ of a discrete Coxeter group $G$ 
inherits as well a central (spinor) extension $\tilde{G}^+$.  
In general, there are central extensions of alternating groups 
different from the spinor ones (see \cite{Maxwell78} for Schur multipliers of alternating subgroups of finite Coxeter groups and for a presentation \`a la Bourbaki). 

\vskip .2cm
We suggest another presentation of the groups $G^+$ and $\tilde{G}^+$ for an arbitrary Coxeter system $(G,S)$. In the Bourbaki presentation the generators are indexed by vertices (except the distinguished one) of the Coxeter graph. In our presentation the generators are indexed by edges of the Coxeter graph.
{}For type A our presentations of the chains of the groups $A_n^+$ and $\tilde{A}_n^+$ are local and stationary. For the chain of the groups $A_n^+$ a similar (but not exactly the same) presentation was given in \cite{VV} and it was asked there how to generalize this presentation to all alternating groups.
{}For types B and D, our presentations 
of the chains of groups 
$B^+_n$, $\tilde{B}_n^+$, $D^+_n$ and $\tilde{D}_n^+$
are local and eventually stationary (that is, stationary for $n\geqslant n_0$ for a certain $n_0$).

\vskip .2cm
The group $A_n^+$ has a presentation referred to as {\it Carmichael presentation} in \cite{CM}. We interpret this presentation geometrically:
here one distinguishes an oriented edge in the Coxeter graph. With this interpretation there are analogues 
of the Carmichael generators for other 
Coxeter systems. In general, these elements generate 
a subgroup of $G^+$. However for type B (with an appropriate choice of a distinguished oriented edge) and type D the whole group $G^+$ is obtained. We give defining relations for these sets of generators of $B_n^+$ and $D^+_n$. 

\vskip .2cm
The paper is organized as follows. In Section \ref{sec-group} we obtain the presentation of the alternating subgroup $G^+$ for a Coxeter system $(G,S)$.
We first deal with irreducible Coxeter systems. Then we generalize the presentation to arbitrary Coxeter systems by introducing a \emph{connected extension} of 
the Coxeter graph.  
Section \ref{sec-spin} is devoted to the presentation of the group $\tilde{G}^+$ for an arbitrary  Coxeter system $(G,S)$.
One can give presentations - with generators related to edges of the Coxeter graph - for all central extensions of the alternating subgroups of finite Coxeter groups. We illustrate it on two examples: $A_5^+$ and $A_6^+$. 
In Section \ref{app} we give results of the Coxeter--Todd algorithm
for three presentations - the Bourbaki one, our presentation which refers to edges and the presentation \`a la Carmichael - of the alternating groups of types A, B and D.

\vskip .2cm
In the text whenever a product $\pi\sigma$ of two permutations appears we assume that $\sigma$ is applied first and then $\pi$; for example,
$(1,2)(2,3)=(1,2,3)$.

\section{{\hspace{-0.55cm}.\hspace{0.55cm}}Alternating subgroups of Coxeter groups}\label{sec-group}

\paragraph{1.} Let $(G,S)$ be a Coxeter system: $S$ is the set of generators, $S=\{ s_0,\dots,s_{n-1}\}$;
the defining relations are encoded by a symmetric matrix $\mathfrak{m}=(m_{ij})_{i,j=0,1,\dots ,n-1}$ with $m_{ii}=1$ and $2\leq m_{ij}\in\mathbb{N}\cup\infty$:
\begin{equation}\label{codecoxet} G:=\left\langle\ s_0,\dots,s_{n-1}\ |\ (s_is_j)^{m_{ij}}=1\ \textrm{for} \ i,j=0,1,\dots ,n-1,\ i\leqslant j\ \right\rangle.\end{equation} 
This presentation of $G$ can be expressed with the help of a {\it Coxeter} graph: its vertices are in one-to-one correspondence with the generators $s_0,\dots,s_{n-1}$ and are indexed by the subscripts $0,1,\dots ,n-1$; 
the vertex $i$ is connected to the vertex $j$ if and only if $m_{ij}\geq3$; the edge between $i$ and $j$ is labeled by the number $m_{ij}$ if $m_{ij}\geq3$
(for brevity, in the sequel we draw a simple line for an edge labeled by $3$ and a double line for an edge labeled by $4$). 
We denote the Coxeter graph corresponding to the Coxeter system $(G,S)$ by ${\cal{G}}$. 

\vskip .1cm
The sign character is the unique homomorphism $\epsilon : G\to \left\{-1,1\right\}$ such that $\epsilon(s_i)=-1$ for $i=0,\dots,n-1$. Its kernel $G^+:=\text{ker}(\epsilon)$ is called the alternating subgroup of $G$. Recall the Bourbaki \cite{Bour} presentation of $G^+$, see \cite{BRR} for a proof.
The alternating group $G^+$ is isomorphic to the group generated by $R_1,\dots,R_{n-1}$ with the defining relations:
\begin{equation}\label{codebourb}\left\{\begin{array}{ll}
R_i^{m_{0i}}=1 & \textrm{for $i=1,\dots,n-1$,}\\[.2em]
(R_i^{-1}R_j)^{m_{ij}}=1 & \textrm{for $i,j=1,\dots,n-1$ such that $i<j$.}\end{array}\right.\end{equation}
The isomorphism with $G^+$ is given by $R_i\mapsto s_0s_i$ for $i=1,\dots,n-1$. 
The Bourbaki presentation of the alternating group $G^+$ depends on the choice of a generator carrying the subscript 0.

\paragraph{2.} In this paragraph, we assume that the Coxeter system $(G,S)$ is irreducible (that is, the Coxeter graph ${\cal{G}}$ is connected). We give a presentation of $G^+$ in terms of generators corresponding to the edges of the Coxeter graph; no vertex in the Coxeter graph is distinguished.
In the next paragraph this presentation will be generalized to arbitrary Coxeter systems.

\vskip .1cm
The presentation uses an orientation - chosen arbitrarily - of edges 
of the Coxeter graph. For concreteness, if there is an edge between $i$ and $j$ with $i<j$,
we orient it from $i$ to $j$. We associate a generator $r_{ij}$ to every oriented edge, that is, to every pair $(i,j)$, $i,j=0,\dots,n-1$, such that $i<j$ and $m_{ij}\neq 2$.
For a generator $r_{ij}$ we denote by $r_{ji}$ the inverse, $r_{ji}:=r_{ij}^{-1}$.

\begin{defi}
\label{def-dis-edge}
Two generators $r_{ij}$ and $r_{lm}$ are said to be not connected if $\left\{i,j\right\}\cap\left\{l,m\right\}=\varnothing$ and there is no edge connecting any of 
the vertices $\left\{i,j\right\}$ with any of the vertices $\left\{l,m\right\}$.
\end{defi}

\begin{prop}\label{prop-code}Let $(G,S)$ be an irreducible Coxeter system with a Coxeter matrix $\mathfrak{m}$. The alternating subgroup $G^+$ of $G$ is isomorphic to the group with  generators $r_{ij}$ and the defining relations
\begin{empheq}[left=\empheqlbrace]{alignat=1}\label{codenotre-a}&(r_{ij})^{m_{ij}}=1 \hspace{2.2cm} \textrm{for all generators $r_{ij}$,}\\[.1em]
\label{codenotre-b}&r_{ii_1}r_{i_1i_2}\dots r_{i_ai}=1 \hspace{1cm} \textrm{for cycles with edges  $(ii_1),(i_1i_2)\dots,(i_ai)$}, 
\\[.1em]
\label{codenotre-c}&(r_{ij}r_{jk})^2=1 \hspace{2cm} \textrm{for $r_{ij},r_{jk}$ such that $i<k$ and $m_{ik}=2$,}\\[.1em]
\label{codenotre-d}&(r_{ij}r_{jk}r_{kl})^2=1 \hspace{1.5cm} \textrm{for $r_{ij},r_{jk},r_{kl}$ such that $i<l$ and $m_{il}=2$,}\\[.1em]
\label{codenotre-e}&r_{ij}r_{lm}=r_{lm}r_{ij} \hspace{1.55cm} \textrm{for $r_{ij}$ and $r_{lm}$ which are not connected.}\end{empheq}
\end{prop}
Let $\mathfrak{c}_1,\dots,\mathfrak{c}_{\mathfrak{l}}$ be a set of generators of the fundamental group of ${\cal{G}}$. In the set of the defining relations 
it is sufficient to impose the relation (\ref{codenotre-b}) for the cycles $\mathfrak{c}_{\mathfrak{a}}$, $\mathfrak{a}=1,\dots,\mathfrak{l}$. If $G$ is finite, its Coxeter graph has no cycles so the relation (\ref{codenotre-b}) is not needed. 

\noindent\textit{Proof of the Proposition.} Let $W^+$ be the group generated by $r_{ij}$ with the defining relations (\ref{codenotre-a})--(\ref{codenotre-e}). Define the map $\phi$ from the set of generators $r_{ij}$ of $W^+$ to $G^+$ by
\begin{equation}\label{real}\phi\ :\ r_{ij}\mapsto \left\{\begin{array}{ll}R_i^{-1}R_j & \textrm{if $i\neq0$,}\\[.2em] 
R_j & \textrm{if $i=0$.}\end{array}\right.\end{equation}
It is straightforward to see that the relations (\ref{codenotre-a})--(\ref{codenotre-d}) are verified by the images of the generators
$r_{ij}$ in $G^+$. As for the last one, the generators $r_{ij}$ and $r_{lm}$ are not connected, so $m_{il}=m_{im}=m_{jl}=m_{jm}=2$. 
We have two possibilities: either $i\neq 0$ and $l\neq 0$ or one of the numbers, $i$ or $l$, is $0$.
\begin{itemize}
\item If $i\neq 0$ and $l\neq 0$ then
\begin{equation}\label{calc}R_i^{-1}R_jR_l^{-1}R_m=R_i^{-1}R_lR_j^{-1}R_m=R_l^{-1}R_iR_j^{-1}R_m=R_l^{-1}R_iR_m^{-1}R_j=R_l^{-1}R_mR_i^{-1}R_l,\end{equation}
where we have repeatedly used the second line in (\ref{codebourb}).
Thus $\phi(r_{ij})\phi(r_{lm})=\phi(r_{lm})\phi(r_{ij})$.
 \item If, say, $l=0$ then $i\neq 0$ and we have $R_i^2=R_j^2=1$ because $m_{0i}=m_{0j}=2$. We calculate
\begin{equation}\label{calc2} R_i^{-1}R_jR_m=R_i^{-1}R_m^{-1}R_j=R_mR_i^{-1}R_j;\end{equation}
we used  $(R_j^{-1}R_m)^2=1$ and $R_j^2=1$ in the first equality; in the second equality we used $(R_i^{-1}R_m)^2=1$ and $R_i^2=1$. Thus $\phi(r_{ij})\phi(r_{lm})=\phi(r_{lm})\phi(r_{ij})$ again.
\end{itemize}
We conclude that $\phi$ extends to a group homomorphism from $W^+$ to $G^+$.

Define now the map $\psi$ from the set of generators of $G^+$ to $W^+$ by
\begin{equation}\label{inv-real}\psi\ :\ R_i\mapsto{\grave{R}}_i:=r_{0i_1}r_{i_1i_2}\dots r_{i_ki}\ \ \textrm{for all $i=1,\dots,n-1$,}\end{equation}
where $(0,i_1,i_2,\dots,i_k,i)$ is an arbitrary path from the vertex $0$ to the vertex $i$ in the Coxeter graph (such path exists because the Coxeter graph is connected). Due to the relation (\ref{codenotre-b}), the element $\grave{R}_i$ does not depend on the chosen path,
thus the map $\psi$ is well-defined. We shall assume that $(0,i_1,i_2,\dots,i_k,i)$ is a shortest (that is, with minimal number of edges) path from $0$ to $i$. 

If $j\neq i$ then
\[\grave{R}_i^{-1}\grave{R}_j=r_{ii_k}\dots r_{i_10}r_{0j_1}\dots r_{j_lj}\ ,\ \text{or, again by (\ref{codenotre-b})},\ \grave{R}_i^{-1}\grave{R}_j=r_{ii_a}\dots r_{i_bj},\] 
where $(i,i_a,\dots,i_b,j)$ is a shortest path from $i$ to $j$. 

We shall verify that $\grave{R}_i^{m_{0i}}=1$ for $i=1,\dots, n-1$, and $(\grave{R}_i^{-1}\grave{R}_j)^{m_{ij}}=1$ for $i,j=1,\dots,n-1$ with $i<j$. 
It suffices to prove that $(r_{ia_1}r_{a_1a_2}\dots r_{a_kj})^{m_{ij}}=1$ for a shortest path $(i,a_1,a_2,\dots, a_k,j)$. 
If $(k_1,k_2,\dots,k_l)$ is a shortest path then there is no edge in the Coxeter graph between the vertices $k_{\mu}$ and $k_{\nu}$ with 
$\vert \mu -\nu\vert >1$. If $m_{ij}\neq 2$ the shortest path is $(i,j)$ and $(r_{ij})^{m_{ij}}=1$. 
If $m_{ij}=2$ and $k=1$ (respectively, $k=2$) then $(r_{ia_1}\dots r_{a_kj})^2=1$ by (\ref{codenotre-c}) (respectively, (\ref{codenotre-d})). If $m_{ij}=2$ and $k\geq 3$ we have:
\[\begin{array}{ll}
(r_{ia_1}r_{a_1a_2}\dots r_{a_kj})^2 & =r_{ia_1}r_{a_1a_2}r_{a_2a_3}\dots r_{a_kj}r_{ia_1}\dots r_{a_kj} =r_{a_2a_3}^{-1}r_{a_1a_2}^{-1}r_{ia_1}^{-1}\dots r_{a_kj}r_{ia_1}\dots r_{a_kj}\\[.4em]
 & =r_{a_1a_2}r_{a_2a_3}\dots r_{a_kj}r_{a_1a_2}\dots r_{a_kj} =(r_{a_1a_2}\dots r_{a_kj})^2.
 \end{array}\]
In the second equality we used that $(r_{ia_1}r_{a_1a_2}r_{a_2a_3})^2=1$ since $m_{ia_3}=2$; in the third equality, we used that $r_{ia_1}$ commutes with $r_{a_3a_4}\dots r_{a_kj}$ and, since $m_{a_1a_3}=2$, that $(r_{a_1a_2}r_{a_2a_3})^2=1$. These properties are implied by the
the minimality of the path $(i,a_1,a_2,\dots, a_k,j)$.
The assertion follows by induction on the length of the path.

We conclude that the map $\psi$ extends to a group homomorphism from $G^+$ to $W^+$.
It is straighforward to see that the morphisms $\phi$ and $\psi$ are mutually inverse.\hfill$\square$

\vskip .2cm\noindent
\textbf{Remark.} Let  $(i,a_1,a_2,\dots,a_k,j)$ be any path in the Coxeter graph. As the proof shows,  the relations (\ref{codenotre-a})--(\ref{codenotre-e}) 
imply $(r_{ia_1}r_{a_1a_2}\dots r_{a_kj})^{m_{ij}}=1$. Therefore, the defining relations (\ref{codenotre-a})--(\ref{codenotre-e}) are equivalent to the following relations:
\[ \left\{\begin{array}{lll}
 (r_{ia_1}r_{a_1a_2}\dots r_{a_kj})^{m_{ij}}=1 && \text{for any path $(i,a_1,a_2,\dots,a_k,j)$ in the Coxeter graph,}\\[0.5em]
 r_{ij}r_{lm}=r_{lm}r_{ij} && \textrm{for
$r_{ij}$ and $r_{lm}$ which are not connected.}
  \end{array}\right.\] 

\paragraph{3.} Now assume that the Coxeter graph ${\cal{G}}$ is not connected.
We add some `virtual" edges labeled by $2$ to make the graph connected. More precisely, let ${\cal{G}}={\cal{G}}_1\sqcup{\cal{G}}_2\sqcup\dots\sqcup{\cal{G}}_m$ be a decomposition of ${\cal{G}}$ into the disjoint union of its connected components. 
{}For $a=1,\dots,m$ choose an arbitrary vertex $i_a$
of ${\cal{G}}_a$; then, for all $l=1,\dots,m-1$, we add a virtual edge between $i_l$ and $i_{l+1}$  and label it by $2$. 
Virtual edges will be drawn with dashed lines. The obtained graph we call a \emph{connected extension} of ${\cal{G}}$. 

A connected extension is not unique. Fix a connected extension ${\cal{G}}^c$ of ${\cal{G}}$ 
and associate a generator $r_{ij}$, $i<j$, to every oriented edge  (virtual or not) of ${\cal{G}}^c$.
The analogue of the Definition \ref{def-dis-edge} includes now virtual edges.
\begin{defi}
\label{def-dis-edge2}
Two generators $r_{ij}$ and $r_{kl}$ are said to be not connected if $\left\{i,j\right\}\cap\left\{l,m\right\}=\varnothing$ and there is no edge or virtual edge connecting any of 
the vertices $\left\{i,j\right\}$ with any of the vertices $\left\{l,m\right\}$.
\end{defi}

With this choice of the set of generators, the formulation of the Proposition \ref{prop-code} stays the same for an arbitrary Coxeter system.
\begin{prop}
\label{prop-code'}
Let $(G,S)$ be a Coxeter system with a Coxeter matrix $\mathfrak{m}$. The alternating subgroup $G^+$ of $G$ is isomorphic to the group generated by $r_{ij}$  with the defining relations (\ref{codenotre-a})--(\ref{codenotre-e}).
\end{prop}
\emph{Proof.}
The proof goes exactly as the proof of the Proposition \ref{prop-code}. The isomorphisms are still given by
\begin{equation}\phi\ :\ r_{ij}\mapsto \left\{\begin{array}{ll}R_i^{-1}R_j & \textrm{if $i\neq0$,}\\[.2em] 
R_j & \textrm{if $i=0$,}\end{array}\right.\end{equation}
and
\begin{equation}\psi\ :\ R_i\mapsto
r_{0i_1}r_{i_1i_2}\dots r_{i_ki}\ \ \textrm{for $i=1,\dots,n-1$,}\end{equation}
where $(0,i_1,i_2,\dots,i_k,i)$ is an arbitrary path from the vertex $0$ to the vertex $i$ in ${\cal{G}}^c$.
\hfill$\square$

\vskip .2cm\noindent 
\textbf{Remark.} Denote by $G({\cal{G}})$ the Coxeter group related to the Coxeter graph ${\cal{G}}$. The defining relations (\ref{codenotre-a})--(\ref{codenotre-e}) imply that generators $r_{ij}$ and $r_{kl}$ commute if the edges $(i,j)$ and $(k,l)$ belong to different connected components of ${\cal{G}}$. 
Indeed, either $r_{ij}$ and $r_{kl}$ are not connected and they commute, or there is a virtual edge between, say, vertices $j$ and $k$; in the latter situation, relations (\ref{codenotre-c})--(\ref{codenotre-d}) gives $$r_{kl}r_{ij}=r_{jk}^{-1}r_{ij}^{-1}r_{kl}^{-1}r_{jk}^{-1}=r_{ij}r_{jk}^2r_{kl},$$ and the assertion follows from $r_{jk}^2=1$.

Let, as above, ${\cal{G}}={\cal{G}}_1\sqcup{\cal{G}}_2\sqcup\dots\sqcup{\cal{G}}_m$; the groups $G({\cal{G}}_a)$, $a=1,\dots ,m$, are naturally subgroups of $G({\cal{G}})$. The virtual edges form a path in ${\cal{G}}^c$ and
one can verify that 
the generators $r_{st}$ corresponding to the virtual edges of ${\cal{G}}^c$ generate a subgroup ${\cal{Y}}$ isomorphic to
$C_2^{m-1}$, where $C_2$ is the cyclic group with 2 elements. 
One can also verify that each $G({\cal{G}}_a)$, $a=1,\dots ,m$, is stable under conjugations by elements of ${\cal{Y}}$.
Thus, each $G({\cal{G}}_a)$, $a=1,\dots ,m$, is normal in $G({\cal{G}})$ and $G({\cal{G}})$ is the semi-direct product ${\cal{Y}}\ltimes\left( G({\cal{G}}_1)\times G({\cal{G}}_2)\times\dots\times G({\cal{G}}_m)
\right)$.

\bigskip\noindent
\textbf{Example.} Consider the Coxeter group generated by $s_1,s_2,s_3,s_4,s_5$ with the Coxeter matrix $\mathfrak{m}=\left(\begin{array}{ccccc}1&4&2&2&2\\4&1&2&2&2\\2&2&1&3&3\\2&2&3&1&3\\2&2&3&3&1\end{array}\right)$.
The Coxeter graph has two components which we connect by adding a virtual edge between vertices $2$ and $3$. 
The figure shows oriented edges and generators corresponding to each oriented edge and virtual edge. With the Definition \ref{def-dis-edge2} there is only one pair of generators which are not connected: $r_{12}$ and $r_{45}$.

\begin{center}
\setlength{\unitlength}{2200sp}
\begingroup\makeatletter\ifx\SetFigFontNFSS\undefined
\gdef\SetFigFontNFSS#1#2#3#4#5{
  \reset@font\fontsize{#1}{#2pt}
  \fontfamily{#3}\fontseries{#4}\fontshape{#5}
  \selectfont}
\fi\endgroup
\begin{picture}(9500,2236)(1066,-4580)
{\thinlines
\put(2989,-3346){\circle*{90}}}{\put(5626,-4201){\circle*{90}}
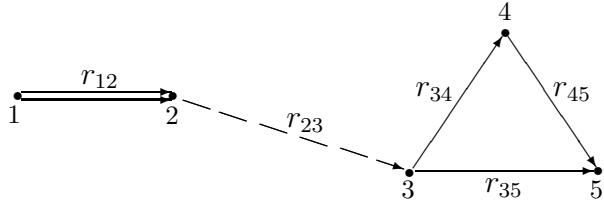
}{\put(6695,-2630){\circle*{90}}}{\put(7730,-4180){\circle*{90}}}{\put(1261,-3346){\circle*{90}}}{\put(3106,-3391){\line( 3,-1){300}}}{\put(3511,-3526){\line( 3,-1){300}}}{\put(3961,-3661){\line( 3,-1){300}}}{\put(4411,-3796){\line( 3,-1){300}}}{\put(4861,-3931){\line( 3,-1){300}}}{\put(5266,-4066){\vector( 3,-1){300}}}{\put(5700,-4180){\vector( 1, 0){1980}}}{\put(5671,-4156){\vector( 2, 3){990}}}{\put(6730,-2670){\vector( 2,-3){980.692}}}{\put(1306,-3300){\vector( 1,0){1665}}}{\put(1306,-3380){\vector(1,0){1665}}}\put(2580,-3655){\makebox(0,0)[lb]{\smash{{\SetFigFontNFSS{12}{14.4}{\rmdefault}{\mddefault}{\updefault}{\small{$2$}}}}}}\put(820,-3640){\makebox(0,0)[lb]{\smash{{\SetFigFontNFSS{12}{14.4}{\rmdefault}{\mddefault}{\updefault}{\small{$1$}}}}}}\put(5210,-4481){\makebox(0,0)[lb]{\smash{{\SetFigFontNFSS{12}{14.4}{\rmdefault}{\mddefault}{\updefault}{\small{$3$}}}}}}\put(7310,-4481){\makebox(0,0)[lb]{\smash{{\SetFigFontNFSS{12}{14.4}{\rmdefault}{\mddefault}{\updefault}{\small{$5$}}}}}}\put(6290,-2491){\makebox(0,0)[lb]{\smash{{\SetFigFontNFSS{12}{14.4}{\rmdefault}{\mddefault}{\updefault}{\small{$4$}}}}}}\put(1636,-3200){\makebox(0,0)[lb]{\smash{{\SetFigFontNFSS{12}{14.4}{\rmdefault}{\mddefault}{\updefault}{$r_{12}$}}}}}\put(3931,-3690){\makebox(0,0)[lb]{\smash{{\SetFigFontNFSS{12}{14.4}{\rmdefault}{\mddefault}{\updefault}{$r_{23}$}}}}}\put(5370,-3346){\makebox(0,0)[lb]{\smash{{\SetFigFontNFSS{12}{14.4}{\rmdefault}{\mddefault}{\updefault}{$r_{34}$}}}}}\put(6901,-3346){\makebox(0,0)[lb]{\smash{{\SetFigFontNFSS{12}{14.4}{\rmdefault}{\mddefault}{\updefault}{$r_{45}$}}}}}\put(6140,-4420){\makebox(0,0)[lb]{\smash{{\SetFigFontNFSS{12}{14.4}{\rmdefault}{\mddefault}{\updefault}{$r_{35}$}}}}}\put(1300,-4900){\begin{pic}\label{Cox-exa}Example of a connected extension of the Coxeter graph
\end{pic}}
\end{picture}
\end{center}

\vskip .4cm
The defining relations for the generators $r_{12}$, $r_{23}$, $r_{34}$, $r_{35}$ and $r_{45}$ are:
\[\left\{\begin{array}{l}(r_{12})^4=1,\  (r_{34})^3=1, \ (r_{35})^3=1,\ (r_{45})^3=1,\ (r_{23})^2=1,\\[.2em]
r_{34}r_{45}r_{53}=1,\\[.2em]
(r_{12}r_{23})^2=1,\ (r_{23}r_{34})^2=1,\ (r_{23}r_{35})^2=1,\\[.2em]
(r_{12}r_{23}r_{34})^2=1,\ (r_{12}r_{23}r_{35})^2=1,\ (r_{23}r_{34}r_{45})^2=1,\ (r_{23}r_{35}r_{54})^2=1,\\[.2em]
r_{12}r_{45}=r_{45}r_{12}.\end{array}\right.\]

\paragraph{4.} We write explicitly the presentation of the Proposition \ref{prop-code} for the alternating groups of type $A$. The Coxeter group $A_n$ is isomorphic to the group of permutations of $n+1$ elements (the symmetric group) and the alternating group $A^+_n$ is isomorphic to the subgroup of even permutations. We use the notation for the generators explained by the following figure (recall that a vertex $i$ corresponds to a generator $s_i$ of the Coxeter group and that a generator of the alternating group is associated to each oriented edge once an orientation is chosen; we have set $r_i:=r_{i-1,i}$ for all $i=1,\dots,n-1$):

\bigskip
\begin{center}
\setlength{\unitlength}{2000sp}
\begingroup\makeatletter\ifx\SetFigFontNFSS\undefined
\gdef\SetFigFontNFSS#1#2#3#4#5{
  \reset@font\fontsize{#1}{#2pt}
  \fontfamily{#3}\fontseries{#4}\fontshape{#5}
  \selectfont}
\fi\endgroup
\begin{picture}(10000,694)(796,-3779)
{\thinlines\put(2989,-3346){\circle*{144}}}{\put(1171,-3346){\circle*{144}}}{\put(10081,-3346){\circle*{144}}}{\put(4861,-3346){\line(1,0){270}}}{\put(5356,-3346){\line(1,0){270}}}{\put(5851,-3346){\line(1,0){270}}}{\put(6346,-3346){\line(1,0){270}}}{\put(6886,-3346){\line(1,0){270}}}{\put(7381,-3346){\line(1,0){270}}}{\put(7876,-3346){\line(1,0){270}}}{\put(2881,-3346){\vector(-1,0){1665}}}{\put(4726,-3346){\vector(-1,0){1665}}}{\put(10036,-3346){\vector(-1,0){1665}}}\put(8656,-3156){\makebox(0,0)[lb]{\smash{{\SetFigFontNFSS{12}{14.4}{\rmdefault}{\mddefault}{\updefault}{$r_1$}}}}}\put(2380,-3706){\makebox(0,0)[lb]{\smash{{\SetFigFontNFSS{12}{14.4}{\rmdefault}{\mddefault}{\updefault}{\small{$n-2$}}}}}}\put(3346,-3156){\makebox(0,0)[lb]{\smash{{\SetFigFontNFSS{12}{14.4}{\rmdefault}{\mddefault}{\updefault}{$r_{n-2}$}}}}}\put(1456,-3156){\makebox(0,0)[lb]{\smash{{\SetFigFontNFSS{12}{14.4}{\rmdefault}{\mddefault}{\updefault}{$r_{n-1}$}}}}}\put(9740,-3706){\makebox(0,0)[lb]{\smash{{\SetFigFontNFSS{12}{14.4}{\rmdefault}{\mddefault}{\updefault}{\small{$0$}}}}}}\put(520,-3706){\makebox(0,0)[lb]{\smash{{\SetFigFontNFSS{12}{14.4}{\rmdefault}{\mddefault}{\updefault}{\small{$n-1$}}}}}}\put(3000,-4300){\begin{pic}\label{Cox-A}Coxeter graph of type A
\end{pic}}
\end{picture}
\end{center}

\vskip.4cm
According to the Proposition \ref{prop-code} the group $A^+_n$ is generated by $r_i$, $i=1,\dots,n-1$, with the defining relations:
\begin{equation}\label{code1}\left\{\begin{array}{ll}
{r_i}^3=1 & \textrm{for $i=1,\dots,n-1$,}\\[.2em]
(r_i r_{i+1})^2=1 & \textrm{for $i=1,\dots,n-2$,}\\[.2em]
(r_i r_{i+1}r_{i+2})^2=1 & \textrm{for $i=1,\dots,n-3$,}\\[.2em]
r_ir_j=r_jr_i & \textrm{for $i,j=1,\dots,n-1$ such that $|i-j|>2$.}
\end{array} \right.\end{equation}
The isomorphism with the group of even permutations of $n+1$ elements is established by $r_i\mapsto (i,i+1,i+2)$ where $(i,i+1,i+2)$ is the cyclic permutation of $i$, $i+1$ and $i+2$. In Section \ref{app} we give a different proof, based on the Coxeter-Todd algorithm, for this presentation of $A^+_n$.

\bigskip
The presentation (\ref{code1}) is a local and stationary, in the sense of \cite{Vershik91}, presentation for the chain of groups $A_n^+$. In \cite{VV} a presentation of the group $A^+_n$ is proved. The generators are $\rho_1,\dots,\rho_{n-1}$ and the defining relations are the same as the defining relations (\ref{code1}) for the $r_i$, except that the third relation in (\ref{code1}) is replaced by $\rho_i\rho_{i+1}^2\rho_{i+2}=\rho_{i+2}\rho_i$ for $i=1,\dots,n-3$. The 
equivalence of the two presentations is given by $\rho_i\mapsto r_i$ for $i=1,\dots,n-1$. Indeed we have
\[(\rho_i\rho_{i+1}\rho_{i+2})^2=\rho_i\rho_{i+1}\cdot\rho_{i+2}\rho_i\cdot\rho_{i+1}\rho_{i+2}=\rho_i\rho_{i+1}\cdot\rho_i\rho_{i+1}^2\rho_{i+2}\cdot\rho_{i+1}\rho_{i+2}=(\rho_i\rho_{i+1})^2(\rho_{i+1}\rho_{i+2})^2=1.\]
Conversely,
\[r_{i+2}r_i=r_{i+1}^{-1}r_{i}^{-1}\cdot r_{i+2}^{-1}r_{i+1}^{-1}=r_ir_{i+1}\cdot r_{i+1}r_{i+2}=r_ir_{i+1}^2r_{i+2}.\]
The Coxeter-Todd algorithm given in Section \ref{app} provides a normal form, different from the normal 
form in \cite{VV}, for the elements of the group $A^+_n$ .

\bigskip\noindent
\textbf{Remark.} The two first relations of (\ref{code1}) imply that:
\[\left\{\begin{array}{ll}r_ir_{i+1}^2r_i=r_ir_{i+1}\cdot r_{i+1}r_i=r_{i+1}^2r_i^2\cdot r_i^2r_{i+1}^2=r_{i+1}^2r_ir_{i+1}^2,\\[0.5em]
r_i^2r_{i+1}r_i^2=r_i^2r_{i+1}^2\cdot r_{i+1}^2r_i^2=r_{i+1}r_i\cdot r_ir_{i+1}=r_{i+1}r_i^2r_{i+1},\end{array}\right.\ \ \ i=1,\dots,n-2\,.\]
Thus the generators $r'_i:=r_i^{-1}$ for $i$ odd and $r'_i:=r_i$ for $i$ even 
verify the Artin relation 
$$r'_ir'_{i+1}r'_i=r'_{i+1}r'_ir'_{i+1}\quad\textrm{for $i=1,\dots,n-2$.}$$

\section{{\hspace{-0.55cm}.\hspace{0.50cm}}Central extensions of the alternating subgroups of Coxeter groups} \label{sec-spin}

\paragraph{1.} 
Let $(G,S)$ be an arbitrary Coxeter system with $S=\{ s_0,\dots,s_{n-1}\}$ 
and the Coxeter matrix $\mathfrak{m}$. 
Let $\tilde{G}$ be the group with the generators $\tilde{s}_0,\dots,\tilde{s}_{n-1}$, $\alpha$ and  the defining relations
\begin{equation}\label{code-ext}\left\{\begin{array}{ll}
(\tilde{s}_i\tilde{s}_j)^{m_{ij}}=1\qquad \textrm{if $m_{ij}$ is odd,}\\[.2em]
(\tilde{s}_i\tilde{s}_j)^{m_{ij}}=\alpha\qquad \textrm{if $m_{ij}$ is even},\\[.2em]
\textrm{$\alpha$ is central and $\alpha^2=1$}.
\end{array}\right.\end{equation}
Let also $\tilde{G}'$ be the group with the generators $\tilde{s}_0',\dots,\tilde{s}_{n-1}'$, $\alpha'$ and the defining relations
\begin{equation}\label{code-ext'}\left\{\begin{array}{ll} 
(\tilde{s}_i'\tilde{s}_j')^{m_{ij}}=\alpha',\\[.2em]
\textrm{$\alpha'$ is central and $\alpha'^2=1$}.
\end{array}\right.\end{equation}
When $G$ is a discrete reflection group, $\tilde{G}$ and $\tilde{G}'$ are its spinor central extensions and they are non-trivial if there are some $i,j$, $i\neq j$, such that 
$m_{ij}$ is even \cite{Mo2}.
In general, the groups $\tilde{G}$ and $\tilde{G}'$ are not isomorphic. 
Note that discrete Coxeter groups admit non-trivial central extensions 
non-isomorphic to (\ref{code-ext}) or (\ref{code-ext'}); see \cite{IY,Y} for the Schur multipliers of the discrete Coxeter groups.

Denote by $\pi\colon \tilde{G}\rightarrow G$ and $\pi'\colon \tilde{G}'\rightarrow G$ the natural projections. Let $\tilde{G}^+:=\pi^{-1}(G^+)$ and
$\tilde{G}'^+:=\pi'^{-1}(G^+)$. 

\vskip .2cm
Our aim is to give presentations for the groups $\tilde{G}^+$ and $\tilde{G}'^+$ in the spirit of the preceding Section. As a by-product, we 
shall see that the groups $\tilde{G}^+$ and $\tilde{G}'^+$ are isomorphic.

\vskip .2cm
We first give a presentation for $\tilde{G}^+$. 
\begin{prop}
\label{prop-code2-bour}
Let $(G,S)$ be a Coxeter system with the Coxeter matrix $\mathfrak{m}$. The group $\tilde{G}^+$ is generated by $\tR_1,\dots,\tR_{n-1}$ and $z$ with the defining relations
\begin{empheq}[left=\empheqlbrace]{alignat=1}
\label{codebour2-a}&(\tR_i)^{m_{0i}}=z^{m_{0i}-1} \hspace{1.8cm} \textrm{for $i=1,\dots,n-1$,}\\[.1em]
\label{codebour2-b}&(\tR_i^{-1}\tR_j)^{m_{ij}}=z^{m_{ij}-1} \hspace{1cm} \textrm{for $i,j=1,\dots,n-1$ such that $i<j$,}\\[.1em]
\label{codebour2-c}&\textrm{$z$ is central and } z^2=1.
\end{empheq}
\end{prop}
\emph{Proof.} This is an adaptation of the proof,  suggested in \cite{Bour},
of the presentation (\ref{codebourb}) for the group $G^+$; we give a sketch. Let $\tilde{W}^+$ be the group generated by the elements $\tR_1,\dots,\tR_{n-1}$ and $z$ with the defining relations (\ref{codebour2-a})--(\ref{codebour2-c}).
Define the map $\phi$ from the set of generators $\{ \tR_1,\dots,\tR_{n-1},z\}$ to $\tilde{G}^+$ by: 
\begin{equation}\label{real2}\phi(\tR_{i})=\tilde{s}_0\tilde{s}_i\ ,\textrm{ $i=1,\dots,n-1$, and } \phi(z)=\alpha.\end{equation}
This map extends to a surjective homomorphism from $\tilde{W}^+$ to $\tilde{G}^+$, which we denote by the same symbol $\phi$. We will show that $\phi$ is an isomorphism.

One checks that the map $\omega$ given by $\omega(\tR_i)=\tR_i^{-1}$, $i=1,\dots,n-1$, and $\omega(z)=z$ extends to an involution of $\tilde{W}^+$, defining thereby the semi-direct product $C_2\ltimes \tilde{W}^+$, where $C_2$ is the cyclic group of order $2$. Let $s$ be the generator of $C_2$. One verifies that the following maps
\[\tilde{s}_0\mapsto s,\ \tilde{s}_i\mapsto s\tR_i,\ i=1,\dots,n-1,\ \alpha\mapsto z,\]
\[s\mapsto\tilde{s}_0,\ \tR_i\mapsto\tilde{s}_0\tilde{s}_i,\ i=1,\dots,n-1,\ z\mapsto\alpha,\]
extend to homomorphisms $\psi_1\ \colon\ \tilde{G}\to C_2\ltimes \tilde{W}^+$ and $\psi_2\ \colon\ C_2\ltimes \tilde{W}^+\to\tilde{G}$ such that $\psi_1\psi_2={\rm Id}_{C_2\ltimes \tilde{W}^+}$ and $\psi_2\psi_1={\rm Id}_{\tilde{G}}$. The restriction of $\psi_1$ to $\tilde{G}^+$ induces the homomorphism inverse to $\phi$. \hfill$\square$

\vskip .2cm
Recall our convention for the orientation of edges of a Coxeter graph or of a connected extension of it: if there is an edge between vertices $i$ and $j$, $i<j$, then it is oriented from $i$ to $j$.
Associate a generator $\tr_{ij}$ of $\tilde{G}^+$ to each generator $r_{ij}$ of $G^+$. Denote by $\tr_{ji}$  the inverse of $\tr_{ij}$, $\tr_{ji}:=\tr_{ij}^{-1}$.  Extend the Definition \ref{def-dis-edge2} to the generators $\tr_{ij}$.

\begin{prop}
\label{prop-code2}
Let $(G,S)$ be a Coxeter system with the Coxeter matrix $\mathfrak{m}$. The group $\tilde{G}^+$ is generated by $\tr_{ij}$ and $z$ with the defining relations
\begin{equation*}
\left\{\begin{array}{ll}(\tr_{ij})^{m_{ij}}=z^{m_{ij}-1} & \textrm{for all generators $\tr_{ij}$,}\\[.5em]
\tr_{ii_1}\tr_{i_1i_2}\dots \tr_{i_ai}=1 &\textrm{for cycles with edges  $(ii_1),(i_1i_2)\dots,(i_ai)$}, \\[.5em]
(\tr_{ij}\tr_{jk})^2=z &\textrm{for
$\tr_{ij},\tr_{jk}$ such that $i<k$ and $m_{ik}=2$,}\\[.5em]
(\tr_{ij}\tr_{jk}\tr_{kl})^2=z &\textrm{for
$\tr_{ij},\tr_{jk},\tr_{kl}$ such that $i<l$ and $m_{il}=2$,}\\[.5em]
\tr_{ij}\tr_{lm}=\tr_{lm}\tr_{ij} &\textrm{for
$\tr_{ij}$ and $\tr_{lm}$ which are not connected.}\\[.5em]
\textrm{$z$ is central and }z^2=1.\end{array}\right.
\end{equation*}
\end{prop}

We omit the proof as it goes along the same lines as the proof of the Proposition \ref{prop-code}. 

\vskip .2cm
We give without details the presentations for  $\tilde{G}'^+$. The generators of the presentation for $\tilde{G}'^+$ analogous to the one given in the Proposition \ref{prop-code2-bour} will be denoted by $\tR'_1,\dots,\tR'_{n-1}$ and $z'$. The defining relations are
$$\left\{\begin{array}{l}(\tR'_i)^{m_{0i}}=z' \hspace{1.8cm} \textrm{for $i=1,\dots,n-1$,}\\[.5em]
(\tR_i'^{-1}\tR'_j)^{m_{ij}}=z' \hspace{1cm} \textrm{for $i,j=1,\dots,n-1$ such that $i<j$,}\\[.5em]
 \textrm{$z'$ is central and }z'^2=1.\end{array}\right. $$
For $\tilde{G}'^+$, the generators of the presentation analogous to the one given in the Proposition \ref{prop-code2} will be denoted by $\tr'_{ij}$ and $z'$ with the same conventions as above. The defining relations are 
\begin{equation*}
\left\{\begin{array}{ll}(\tr'_{ij})^{m_{ij}}=z' & \textrm{for all generators $\tr'_{ij}$,}\\[.5em]
\tr'_{ii_1}\tr'_{i_1i_2}\dots \tr'_{i_ai}=z'^{a+1} &\textrm{for cycles with edges  $(ii_1),(i_1i_2)\dots,(i_ai)$}, 
\\[.5em]
(\tr'_{ij}\tr'_{jk})^2=z' &\textrm{for
$\tr'_{ij},\tr'_{jk}$ such that $i<k$ and $m_{ik}=2$,}\\[.5em]
(\tr'_{ij}\tr'_{jk}\tr'_{kl})^2=z' &\textrm{for
$\tr'_{ij},\tr'_{jk},\tr'_{kl}$ such that $i<l$ and $m_{il}=2$,}\\[.5em]
\tr'_{ij}\tr'_{lm}=\tr'_{lm}\tr'_{ij} &\textrm{for
$\tr'_{ij}$ and $\tr'_{lm}$ which are not connected.}\\[.5em]
 \textrm{$z'$ is central and }z'^2=1.\end{array}\right.
\end{equation*}

\begin{prop}
\label{prop-iso} The extensions $\tilde{G}^+$ and $\tilde{G}'^+$ of the group $G^+$ are isomorphic.\end{prop}
\emph{Proof.} The homomorphism  $\tilde{G}^+\to\tilde{G}'^+$ defined on generators by $\tr_{ij}\mapsto z'\tr_{ij}'$ and $z\mapsto z'$ establishes the required isomorphism. \hfill$\square$

\paragraph{2.} Fix, for type A, the numbering of the vertices and orientation of the edges as shown in Fig.~\ref{Cox-A} and set $\tr_i:=\tr_{i-1,i}$, $i=1,\dots,n-1$. The group  $\tilde{A}_n^+$ is the spinor cover of the alternating group $A_n^+$. The presentation for $\tilde{A}_n^+$ from the Proposition \ref{prop-code2} reads
\begin{equation}\label{code2}\left\{\begin{array}{ll}
\tilde{r}_i^3=1,\\[.1em]
(\tilde{r}_i \tilde{r}_{i+1})^2=z,\\[.1em]
(\tilde{r}_i \tilde{r}_{i+1}\tilde{r}_{i+2})^2=z,\\[.1em]
\tilde{r}_i\tilde{r}_j=\tilde{r}_j\tilde{r}_i & \textrm{if $\mid i-j\mid >2$,}\\[.1em]
\textrm{$z$ is central and $z^2=1$.}
\end{array} \right.\end{equation}
This presentation equips the chain of the spinor extensions of the alternating groups of type A with a structure of a local and stationary tower. 

\paragraph{3.} The Schur multipliers for the alternating groups were calculated by Schur in \cite{S}; Maxwell \cite{Maxwell78} generalized this result to  alternating subgroups of finite Coxeter groups and gave a presentation, in the spirit of (\ref{codebourb}), of the corresponding central extensions.     
It is straightforward - although lengthy - to transform this presentation (for all the central extensions of the alternating groups of all finite Coxeter groups) 
to the presentation which uses generators related to the oriented edges of the Coxeter graph. 

As an example we give, omitting details, the presentations for type A.  
The Schur multiplier of $A^+_n$ is $C_2$ if $n\geqslant 3, n\neq 5,6$, and $C_2\times C_3$ if $n=5,6$; here $C_m$ is the cyclic group with $m$ elements. We describe the central extensions of the groups $A^+_5$ and $A^+_6$ with the kernel $C_2\times C_3$; these extensions are
universal central extensions since each of the groups $A^+_5$ and $A^+_6$ is perfect, that is, it coincides with its commutator subgroup. 
As above, we associate a generator $\tilde{r}_i$, $i=1,\dots ,n-1$, to each oriented edge of the Coxeter graph, see Fig. \ref{Cox-A}. 

\begin{itemize}
\item The universal central extension of $A_5^+$ is generated by $\tilde{r}_1$, $\tilde{r}_2$, $\tilde{r}_3$, $\tilde{r}_4$, $z$ and $\zeta$ with the 
defining relations:
\begin{equation}\label{code2b}\left\{\begin{array}{ll}
\tilde{r}_i^3=1\ \ \textrm{for $i=1,\dots,4$,}\\[.2em]
(\tilde{r}_1 \tilde{r}_2)^2=z,\ \ (\tilde{r}_2\tilde{r}_3)^2=z\zeta,\ \ (\tilde{r}_3\tilde{r}_4)^2=z,\\[.2em]
(\tilde{r}_1 \tilde{r}_2\tilde{r}_3)^2=z,\ \ (\tilde{r}_2 \tilde{r}_3\tilde{r}_4)^2=z,\\[.2em]
\tilde{r}_1\tilde{r}_4=\zeta^2\tilde{r}_4\tilde{r}_1,\\[.2em]
\textrm{$z^2=1$ and $z$ is central,}\ \ \textrm{$\zeta^3=1$ and $\zeta$ is central.}
\end{array} \right.\end{equation}

\item The universal central extension of $A_6^+$ is generated by $\tilde{r}_1$, $\tilde{r}_2$, $\tilde{r}_3$, $\tilde{r}_4$, $\tilde{r}_5$, $z$ and $\zeta$ with the defining relations:
\begin{equation}\label{code2c}\left\{\begin{array}{l}
\tilde{r}_i^3=1\ \ \textrm{for $i=1,\dots,5$,}\\[.2em]
(\tilde{r}_i \tilde{r}_{i+1})^2=z\ \ \textrm{for $i=1,\dots,4$,}\\[.2em]
(\tilde{r}_1 \tilde{r}_2\tilde{r}_3)^2=z,\ \ (\tilde{r}_2\tilde{r}_3\tilde{r}_4)^2=z\zeta,\ \ (\tilde{r}_3\tilde{r}_4\tilde{r}_5)^2=z,\\[.2em]
\tilde{r}_1\tilde{r}_4=\zeta\tilde{r}_4\tilde{r}_1,\ \ \tilde{r}_2\tilde{r}_5=\zeta\tilde{r}_5\tilde{r}_2,\ \ \tilde{r}_1\tilde{r}_5=\zeta^2\tilde{r}_5\tilde{r}_1,\\[.2em]
\textrm{$z^2=1$ and $z$ is central,}\ \ \textrm{$\zeta^3=1$ and $\zeta$ is central.} 
\end{array} \right.\end{equation}\end{itemize}

\section{{\hspace{-0.55cm}.\hspace{0.55cm}}\hspace{-.42cm} Coxeter--Todd algorithms and normal forms for alternating groups of types A, B and D} \label{app}

Let $G$ be a finite group with a given presentation by generators $g_1,\dots,g_m$ and the set of defining relations ${\cal{R}}$.
Let $I$ be a subset of $\{1,\dots,m\}$ and $H$ the subgroup of $G$ generated by $g_a$, $a\in I$. 
The  Coxeter--Todd algorithm for the pair $(G,H)$ constructs the set of the left cosets of $H$ in $G$ and the action of the generators on this set \cite{CT}. The result of a  Coxeter--Todd algorithm is a figure whose vertices are labeled by left cosets and the arrows stand for the action of the generators. The algorithm starts with the left coset $H$; only the generators $g_a$, $a\notin I$, may act non-trivially on this coset and give new vertices. At each step we analyze, using the relations from ${\cal{R}}$, the action of the generators on vertices and draw new vertices or identify existing vertices. 
The algorithm is finished when we know the action of all generators on every coset in the figure.

\vskip .2cm
The  Coxeter--Todd algorithm for $(G,H)$ lists the left cosets and thus provides a normal form for elements of $G$ with respect to $H$.
The algorithm implies an upper bound for the cardinality of $G$. 
Namely, let ${\cal{H}}$ be the abstract group with the generators $g_a$, $a\in I$; the set of defining relations for ${\cal{H}}$ is 
the subset of ${\cal{R}}$ consisting of those relations which involve only the generators $g_a$,  $a\in I$. 
There is a natural surjection  ${\cal{H}}\rightarrow H$ therefore the cardinality of $G$ is less or equal than the numbers of vertices in the figure times 
the cardinality of ${\cal{H}}$. For the chain $\{ 1\}=G_0\subset G_1\subset\dots \subset G_n\subset\dots$ of groups, we construct recursively the global normal form for elements of any $G_n$ using the normal form for elements of $G_k$ with respect to $G_{k-1}$, 
$k=1,2,\dots ,n$.

\vskip .2cm 
Here we exhibit results of the Coxeter--Todd algorithm for the chains of the alternating groups of types A, B and D with the presentation of the Proposition \ref{prop-code}. As a matter of comparison we also give in each case Coxeter--Todd figures for the alternating groups with the presentation (\ref{codebourb}), for the alternating groups with a presentation \`a la Carmichael 
and for the Coxeter groups with the standard presentation (\ref{codecoxet}).
Each time the normal form follows straightforwardly from the  Coxeter--Todd figure; we 
illustrate it on the example of the chain of the alternating groups of type A.

The subgroup chosen for realizing the Coxeter--Todd algorithm is always denoted by $H$. In a Coxeter--Todd figure, when the action of a generator $g$ on a coset 
$Y$ is not specified it means that $gY=Y$.
For a generator of order 2 an unoriented edge represents a pair of oppositely oriented edges. 
The action of a generator of order $3$ (respectively, $4$) is often represented by an oriented triangle (respectively, quadrilateral); the generator is written inside it.
We indicate the coset (in the form $uH$ with $u$ a word in the generators) at each vertex of the  Coxeter--Todd figure only for type A; in general, it can be easily found following the edges from the vertex $H$ in the figures.

\subsection{Type A}

The vertices and oriented edges of the Coxeter graph of type A are labeled as on Fig. \ref{Cox-A}, Section \ref{sec-group}. 
The Coxeter group $A_n$ is generated by $s_0$, $\dots$, $s_{n-1}$ with the defining relations
\begin{equation}\label{Sn}
\left\{ \begin{array}{ll}
s_i^2=1 & \text{for $i=0,\dots,n-1$,}\\
s_is_{i+1}s_i=s_{i+1}s_is_{i+1} & \text{for $i=0,\dots,n-2$},\\
s_is_j=s_js_i & \textrm{for $i,j=0,\dots,n-1$ such that $|i-j|>1$.}
\end{array}\right.
\end{equation}
The group $A_n$ is isomorphic to the group of permutations of the set $\{ 1,2,\dots ,n+1\}$. 
The isomorphism is given by $s_i\mapsto(i+1,i+2)$, $i=0,\dots,n-1$.

\vskip .2cm
Let $H$ be the subgroup generated by $s_0$, $\dots$, $s_{n-2}$; here is the  Coxeter--Todd figure for $(A_n,H)$:

\begin{center}
\setlength{\unitlength}{2100sp}
\begingroup\makeatletter\ifx\SetFigFont\undefined
\gdef\SetFigFont#1#2#3#4#5{
  \reset@font\fontsize{#1}{#2pt}
  \fontfamily{#3}\fontseries{#4}\fontshape{#5}
  \selectfont}
\fi\endgroup
\begin{picture}(10800,1400)(616,-3850)
{\thinlines
\put(2989,-3046){\circle*{144}}}{\put(1171,-3046){\circle*{144}}}{\put(10081,-3046){\circle*{144}}}{\put(4861,-3046){\line(1,0){270}}}{\put(5356,-3046){\line(1,0){270}}}{\put(5851,-3046){\line(1,0){270}}}{\put(6346,-3046){\line(1,0){270}}}{\put(6886,-3046){\line(1,0){270}}}{\put(7381,-3046){\line(1,0){270}}}{\put(7876,-3046){\line(1,0){270}}}{\put(1216,-3046){\line(1,0){1665}}}{\put(3061,-3046){\line(1,0){1665}}}{\put(8371,-3046){\line(1,0){1665}}}\put(3381,-2956){\makebox(0,0)[lb]{\smash{{\SetFigFont{12}{14.4}{\rmdefault}{\mddefault}{\updefault}{$s_{n-2}$}}}}}\put(8556,-2956){\makebox(0,0)[lb]{\smash{{\SetFigFont{12}{14.4}{\rmdefault}{\mddefault}{\updefault}{$s_0$}}}}}\put(9276,-3456){\makebox(0,0)[lb]{\smash{{\SetFigFont{12}{14.4}{\rmdefault}{\mddefault}{\updefault}{$s_0\dots s_{n-1}H$}}}}}\put(636,-3456){\makebox(0,0)[lb]{\smash{{\SetFigFont{12}{14.4}{\rmdefault}{\mddefault}{\updefault}{$H$}}}}}\put(2346,-3456){\makebox(0,0)[lb]{\smash{{\SetFigFont{12}{14.4}{\rmdefault}{\mddefault}{\updefault}{$s_{n-1}H$}}}}}\put(1446,-2956){\makebox(0,0)[lb]{\smash{{\SetFigFont{12}{14.4}{\rmdefault}{\mddefault}{\updefault}{$s_{n-1}$}}}}}
\put(500,-4300){\begin{pic}\label{CoxToddA} Coxeter--Todd figure for $(A_n,H)$ for the presentation (\ref{Sn})
\end{pic}}
\end{picture}
\end{center}

\bigskip
We give Coxeter--Todd figures for three presentations of the alternating group $A^+_n$. In each situation, $H$ is the subgroup generated by the first $n-2$  generators. The first two presentations can be found in \cite{CM}. The second one is the presentation (\ref{codebourb}). 
The third one is the presentation (\ref{code1}). 
To illustrate the usefulness of the  Coxeter--Todd algorithm we 
reestablish that these are indeed 
presentations of the group $A^+_n$ and find three normal forms for elements of $A^+_n$.

\paragraph{1.} 
In Carmichael presentation \cite{CM} the generators are $a_1$, $\dots$, $a_{n-1}$ with the defining relations
\begin{equation}
\label{codeCar} \left\{\begin{array}{ll}
a_i^3=1 & \text{for $i=1,\dots,n-1$,}\\
(a_i a_j)^2=1 & \textrm{for $i,j=1,\dots,n-1$ such that $i<j$.}\end{array} \right.
\end{equation}

\vskip 1cm
\setlength{\unitlength}{2500sp}
\begingroup\makeatletter\ifx\SetFigFont\undefined
\gdef\SetFigFont#1#2#3#4#5{
  \reset@font\fontsize{#1}{#2pt}
  \fontfamily{#3}\fontseries{#4}\fontshape{#5}
  \selectfont}
\fi\endgroup
\begin{picture}(1000,4920)(1000,-7000)
{\thinlines \put(5221,-1951){\circle*{128}}}{\put(4501,-3706){\circle*{128}}}{\put(8731,-4291){\circle*{128}}}{\put(7921,-2581){\circle*{128}}}{\put(6256,-2581){\circle*{128}}}{\put(6227,-6079){\circle*{128}}}{\put(8641,-4336){\vector(-4,-3){2340}}}{\put(7876,-2716){\line(-1,-2){360}}\put(7516,-3436){\line(0,1){45}}}{\put(5221,-2041){\line(1,-4){270}}}{\put(7516,-3526){\vector(-1,-2){1260}}}{\put(6256,-2581){\vector(-3,-2){1744.615}}}{\put(5491,-3121){\vector(1,-4){720}}}{\put(4546,-3706){\vector(2,-3){1530}}}{\put(6301,-6001){\vector(0,1){3240}}}{\put(6256,-5956){\vector(0,1){3195}}}{\put(6211,-5956){\vector(0,1){3240}}}{\put(6166,-5956){\vector(0,1){3240}}}{\put(6301,-2671){\vector(3,-2){2295}}}{\put(6166,-2536){\vector(-3,2){872.308}}}{\put(6256,-2581){\vector(1,0){1600}}}\put(6661,-2761){\makebox(0,0)[lb]{\smash{{\SetFigFont{12}{14.4}{\rmdefault}{\mddefault}{\updefault}{$a_{n-2}$}}}}}\put(6760,-4786){\rotatebox{60.0}{\makebox(0,0)[lb]{\smash{{\SetFigFont{12}{14.4}{\rmdefault}{\mddefault}{\updefault}{$a_{n-2}$}}}}}}\put(7066,-5400){\rotatebox{40.0}{\makebox(0,0)[lb]{\smash{{\SetFigFont{12}{14.4}{\rmdefault}{\mddefault}{\updefault}{$a_{n-1}$}}}}}}\put(7606,-3800){\rotatebox{330.0}{\makebox(0,0)[lb]{\smash{{\SetFigFont{12}{14.4}{\rmdefault}{\mddefault}{\updefault}{$a_{n-1}$}}}}}}\put(4951,-4651){\rotatebox{310.0}{\makebox(0,0)[lb]{\smash{{\SetFigFont{12}{14.4}{\rmdefault}{\mddefault}{\updefault}{$a_1$}}}}}}\put(4951,-4651){\rotatebox{310.0}{\makebox(0,0)[lb]{\smash{{\SetFigFont{12}{14.4}{\rmdefault}{\mddefault}{\updefault}{$a_1$}}}}}}\put(5536,-3886){\rotatebox{280.0}{\makebox(0,0)[lb]{\smash{{\SetFigFont{12}{14.4}{\rmdefault}{\mddefault}{\updefault}{$a_2$}}}}}}\put(5300,-2200){\rotatebox{330.0}{\makebox(0,0)[lb]{\smash{{\SetFigFont{12}{14.4}{\rmdefault}{\mddefault}{\updefault}{$a_2$}}}}}}\put(8530,-4400){\makebox(0,0)[lb]{\smash{{\SetFigFont{12}{14.4}{\rmdefault}{\mddefault}{\updefault}{$H$}}}}}\put(7750,-2520){\makebox(0,0)[lb]{\smash{{\SetFigFont{12}{14.4}{\rmdefault}{\mddefault}{\updefault}{$a_{n-2}a_{n-1}^2H$}}}}}\put(3700,-1950){\makebox(0,0)[lb]{\smash{{\SetFigFont{12}{14.4}{\rmdefault}{\mddefault}{\updefault}{$a_2a_{n-1}^2H$}}}}}\put(5850,-2401){\makebox(0,0)[lb]{\smash{{\SetFigFont{12}{14.4}{\rmdefault}{\mddefault}{\updefault}{$a_{n-1}^2H$}}}}}\put(2970,-3700){\makebox(0,0)[lb]{\smash{{\SetFigFont{12}{14.4}{\rmdefault}{\mddefault}{\updefault}{$a_1a_{n-1}^2H$}}}}}\put(5750,-6400){\makebox(0,0)[lb]{\smash{{\SetFigFont{12}{14.4}{\rmdefault}{\mddefault}{\updefault}{$a_{n-1}H$}}}}}\put(5600,-1950){\rotatebox{350.0}{\makebox(0,0)[lb]{\smash{{\SetFigFont{12}{14.4}{\rmdefault}{\mddefault}{\updefault}{$\dots$}}}}}}\put(7111,-2250){\rotatebox{340.0}{\makebox(0,0)[lb]{\smash{{\SetFigFont{12}{14.4}{\rmdefault}{\mddefault}{\updefault}{$\dots$}}}}}}
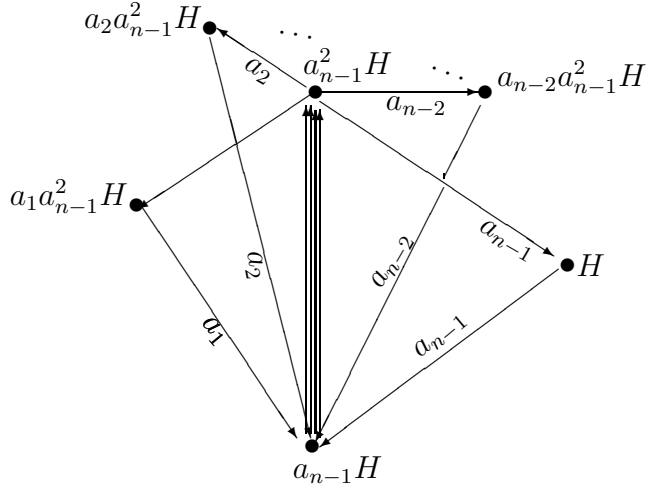
\put(1500,-7050){\begin{pic}\label{CoxToddaltA1} Coxeter--Todd figure for $(A^+_n,H)$ for Carmichael presentation (\ref{codeCar})
\end{pic}}
\end{picture}

\vskip .2cm
We geometrically interpret the Carmichael presentation as follows: fix the oriented edge $(0,1)$ in the Coxeter graph of type A and take $a_1:=s_0s_1$. The other generators are obtained by consecutive conjugations of $a_1$ by the Coxeter generators, namely $a_i:=s_ia_{i-1}s_i$, $i=2,\dots,n-1$. The set $\{a_1,\dots,a_{n-1}\}$ is a generating set for the alternating group $A^+_n$ with the defining relations (\ref{codeCar}).

\paragraph{2.}
Moore presentation \cite{CM} is the presentation (\ref{codebourb}). 
The generators are $R_1$, $\dots$, $R_{n-1}$ with the defining relations
\begin{equation}
\label{codeMoo} \left\{\begin{array}{ll}
R_1^3=1, &\\[.2em]
R_i^2=1 & \textrm{for $i=2\dots,n-1$,}\\[.2em]
(R_i^{-1} R_{i+1})^3=1 & \text{for $i=1,\dots,n-2$,}\\[.2em]
(R_i^{-1} R_j)^2=1 & \textrm{for $i,j=1,\dots,n-1$ such that $|i-j|>1$.}\end{array} \right.
\end{equation}

\bigskip
\setlength{\unitlength}{2000sp}
\begingroup\makeatletter\ifx\SetFigFont\undefined
\gdef\SetFigFont#1#2#3#4#5{
  \reset@font\fontsize{#1}{#2pt}
  \fontfamily{#3}\fontseries{#4}\fontshape{#5}
  \selectfont}
\fi\endgroup
\begin{picture}(10695,3439)(0,-5500)
{\thinlines
\put(2989,-3346){\circle*{144}}}{\put(1171,-3346){\circle*{144}}}{\put(10081,-3346){\circle*{144}}}{\put(11283,-4670){\circle*{144}}}{\put(11316,-2110){\circle*{144}}}{\put(4861,-3346){\line(1,0){270}}}{\put(5356,-3346){\line(1,0){270}}}{\put(5851,-3346){\line(1,0){270}}}{\put(6346,-3346){\line(1,0){270}}}{\put(6886,-3346){\line(1,0){270}}}{\put(7381,-3346){\line(1,0){270}}}{\put(7876,-3346){\line(1,0){270}}}{\put(10126,-3301){\vector(1,1){1170}}}{\put(11296,-2221){\vector(0,-1){2430}}}{\put(11386,-2221){\line(0,-1){2430}}}{\put(11701,-2221){\line(0,-1){2430}}}{\put(11295,-4605){\vector(-1,1){1170}}}{\put(1216,-3346){\line(1,0){1665}}}{\put(3061,-3346){\line(1,0){1665}}}{\put(8371,-3346){\line(1,0){1665}}}\put(8656,-3256){\makebox(0,0)[lb]{\smash{{\SetFigFont{12}{14.4}{\rmdefault}{\mddefault}{\updefault}{\small{$R_2$}}}}}}\put(680,-3706){\makebox(0,0)[lb]{\smash{{\SetFigFont{12}{14.4}{\rmdefault}{\mddefault}{\updefault}{\small{$H$}}}}}}\put(2350,-3706){\makebox(0,0)[lb]{\smash{{\SetFigFont{12}{14.4}{\rmdefault}{\mddefault}{\updefault}{\small{$R_{n-1}H$}}}}}}\put(10250,-3346){\makebox(0,0)[lb]{\smash{{\SetFigFont{12}{14.4}{\rmdefault}{\mddefault}{\updefault}{\small{$R_1$}}}}}}\put(10900,-1906){\makebox(0,0)[lb]{\smash{{\SetFigFont{12}{14.4}{\rmdefault}{\mddefault}{\updefault}{\small{$R_1R_2\dots R_{n-1}H$}}}}}}\put(10750,-5101){\makebox(0,0)[lb]{\smash{{\SetFigFont{12}{14.4}{\rmdefault}{\mddefault}{\updefault}{\small{$R_1^2R_2\dots R_{n-2}H$}}}}}}\put(11400,-3391){\makebox(0,0)[lb]{\smash{{\SetFigFont{12}{14.4}{\rmdefault}{\mddefault}{\updefault}{\small{$\,R_2,\dots,R_{n-1}$}}}}}}\put(1546,-3256){\makebox(0,0)[lb]{\smash{{\SetFigFont{12}{14.4}{\rmdefault}{\mddefault}{\updefault}{\small{$R_{n-1}$}}}}}}\put(3481,-3256){\makebox(0,0)[lb]{\smash{{\SetFigFont{12}{14.4}{\rmdefault}{\mddefault}{\updefault}{\small{$R_{n-2}$}}}}}}\put(7900,-3796){\makebox(0,0)[lb]{\smash{{\SetFigFont{12}{14.4}{\rmdefault}{\mddefault}{\updefault}{\small{$R_2\dots R_{n-1}H$}}}}}}\put(11086,-3256){\makebox(0,0)[lb]{\smash{{\SetFigFont{12}{14.4}{\rmdefault}{\mddefault}{\updefault}{\small{$...$}}}}}}
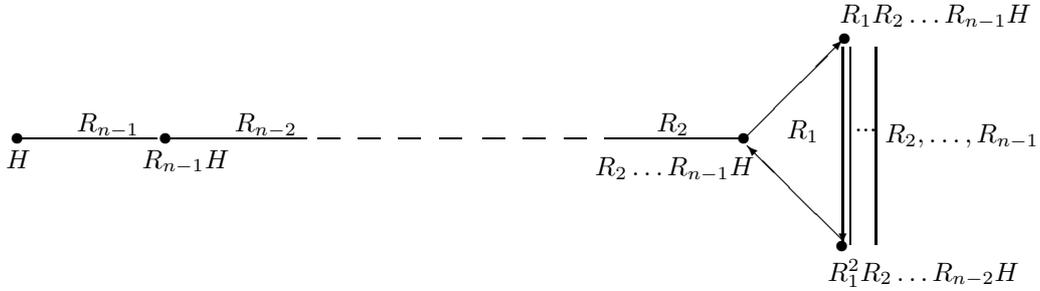
\put(1200,-5750){\begin{pic}\label{CoxToddaltA2} Coxeter--Todd figure for $(A^+_n,H)$ for Moore presentation (\ref{codeMoo})
\end{pic}}
\end{picture}

\paragraph{3.} 
In the presentation (\ref{code1}) of $A^+_n$ the generators are $r_1$, $\dots$, $r_{n-1}$ with the defining relations
\begin{equation}
\label{codeAn} \left\{\begin{array}{ll}
r_i^3=1 & \text{for $i=1,\dots,n-1$,}\\[.2em]
(r_i r_{i+1})^2=1 & \text{for $i=1,\dots,n-2$,}\\[.2em]
(r_i r_{i+1}r_{i+2})^2=1 & \text{for $i=1,\dots,n-3$,}\\[.2em]
r_ir_j=r_jr_i & \textrm{for $i,j=1,\dots,n-1$ such that $|i-j|>2$.}
\end{array} \right.
\end{equation}

\bigskip
\setlength{\unitlength}{2500sp}

\begingroup\makeatletter\ifx\SetFigFontNFSS\undefined
\gdef\SetFigFontNFSS#1#2#3#4#5{
  \reset@font\fontsize{#1}{#2pt}
  \fontfamily{#3}\fontseries{#4}\fontshape{#5}
  \selectfont}
\fi\endgroup
\begin{picture}(10175,3629)(931,-6579)
{\thinlines
\put(2989,-3346){\circle*{144}}}{\put(1171,-3346){\circle*{144}}}{\put(2071,-5146){\circle*{144}}}{\put(4726,-3346){\circle*{144}}}{\put(3916,-5146){\circle*{144}}}{\put(8371,-5146){\circle*{144}}}{\put(9226,-3346){\circle*{144}}}{\put(11026,-3346){\circle*{144}}}{\put(10216,-5191){\circle*{144}}}{\put(3862,-5128){\vector(-1,2){855}}}{\put(3025,-3418){\vector(-1,-2){855}}}{\put(4753,-3382){\vector(-1,-2){855}}}{\put(5131,-3346){\line(1,0){270}}}{\put(5626,-3346){\line(1,0){270}}}{\put(6121,-3346){\line(1,0){270}}}{\put(6661,-3346){\line(1,0){270}}}{\put(7156,-3346){\line(1,0){270}}}{\put(7651,-3346){\line(1,0){270}}}{\put(10090,-5119){\vector(-1,2){855}}}{\put(8191,-3346){\line(1,0){270}}}{\put(8641,-3346){\line(1,0){270}}}{\put(11071,-3436){\vector(-1,-2){855}}}{\put(4276,-5191){\line(1,0){270}}}{\put(4816,-5191){\line(1,0){270}}}{\put(5356,-5191){\line(1,0){270}}}{\put(5896,-5191){\line(1,0){270}}}{\put(6436,-5191){\line(1,0){270}}}{\put(6976,-5191){\line(1,0){270}}}{\put(7516,-5191){\line(1,0){270}}}{\put(7966,-5191){\line(1,0){270}}}{\put(9244,-3400){\vector(-1,-2){855}}}{\put(2035,-5074){\vector(-1,2){855}}}{\put(2971,-3346){\vector(-1,-2){855}}}{\put(10153,-5065){\vector(-1,2){855}}}{\put(3943,-5110){\vector(-1,2){855}}}{\put(2161,-5191){\vector(1,0){1665}}}{\put(9316,-3346){\vector(1,0){1665}}}{\put(8506,-5191){\vector(1,0){1665}}}{\put(1261,-3346){\vector(1,0){1665}}}{\put(3016,-3346){\vector(1,0){1665}}}\put(746,-3166){\makebox(0,0)[lb]{\smash{{\SetFigFontNFSS{12}{14.4}{\rmdefault}{\mddefault}{\updefault}{$H$}}}}}\put(3221,-5506){\makebox(0,0)[lb]{\smash{{\SetFigFontNFSS{12}{14.4}{\rmdefault}{\mddefault}{\updefault}{$r_{n-2}r_{n-1}^2H$}}}}}\put(1421,-5506){\makebox(0,0)[lb]{\smash{{\SetFigFontNFSS{12}{14.4}{\rmdefault}{\mddefault}{\updefault}{$r_{n-1}^2H$}}}}}\put(2321,-3121){\makebox(0,0)[lb]{\smash{{\SetFigFontNFSS{12}{14.4}{\rmdefault}{\mddefault}{\updefault}{$r_{n-1}H$}}}}}\put(10321,-3121){\makebox(0,0)[lb]{\smash{{\SetFigFontNFSS{12}{14.4}{\rmdefault}{\mddefault}{\updefault}{$r_1r_3\dots r_{n-3}r_{n-1}H$}}}}}\put(3931,-3121){\makebox(0,0)[lb]{\smash{{\SetFigFontNFSS{12}{14.4}{\rmdefault}{\mddefault}{\updefault}{$r_{n-3}r_{n-1}H$}}}}}\put(6901,-5506){\makebox(0,0)[lb]{\smash{{\SetFigFontNFSS{12}{14.4}{\rmdefault}{\mddefault}{\updefault}{$r_4\dots r_{n-2}r_{n-1}^2H$}}}}}\put(9556,-5506){\makebox(0,0)[lb]{\smash{{\SetFigFontNFSS{12}{14.4}{\rmdefault}{\mddefault}{\updefault}{$r_2r_4\dots r_{n-2}r_{n-1}^2H$}}}}}\put(7891,-3121){\makebox(0,0)[lb]{\smash{{\SetFigFontNFSS{12}{14.4}{\rmdefault}{\mddefault}{\updefault}{$r_3\dots r_{n-3}r_{n-1}H$}}}}}\put(8801,-4606){\makebox(0,0)[lb]{\smash{{\SetFigFontNFSS{12}{14.4}{\rmdefault}{\mddefault}{\updefault}{$r_2$}}}}}\put(2456,-4651){\makebox(0,0)[lb]{\smash{{\SetFigFontNFSS{12}{14.4}{\rmdefault}{\mddefault}{\updefault}{$r_{n-2}$}}}}}\put(3356,-3976){\makebox(0,0)[lb]{\smash{{\SetFigFontNFSS{12}{14.4}{\rmdefault}{\mddefault}{\updefault}{$r_{n-3}$}}}}}\put(1511,-3976){\makebox(0,0)[lb]{\smash{{\SetFigFontNFSS{12}{14.4}{\rmdefault}{\mddefault}{\updefault}{$r_{n-1}$}}}}}\put(9791,-4021){\makebox(0,0)[lb]{\smash{{\SetFigFontNFSS{12}{14.4}{\rmdefault}{\mddefault}{\updefault}{$r_1$}}}}}
\put(1650,-6200){\begin{pic}\label{CoxToddaltA3a} Coxeter--Todd figure for $(A^+_n,H)$ for the presentation (\ref{codeAn}), $n$ even
\end{pic}}
\end{picture}

\bigskip
\setlength{\unitlength}{2500sp}
\begingroup\makeatletter\ifx\SetFigFontNFSS\undefined
\gdef\SetFigFontNFSS#1#2#3#4#5{
  \reset@font\fontsize{#1}{#2pt}
  \fontfamily{#3}\fontseries{#4}\fontshape{#5}
  \selectfont}
\fi\endgroup
\begin{picture}(10040,3629)(931,-6579)
{\thinlines\put(2989,-3346){\circle*{144}}}{\put(1171,-3346){\circle*{144}}}{\put(2071,-5146){\circle*{144}}}{\put(4726,-3346){\circle*{144}}}{\put(3916,-5146){\circle*{144}}}{\put(8191,-3346){\circle*{144}}}{\put(9991,-3346){\circle*{144}}}{\put(9091,-5191){\circle*{144}}}{\put(10891,-5191){\circle*{144}}}{\put(3862,-5128){\vector(-1,2){855}}}{\put(3025,-3418){\vector(-1,-2){855}}}{\put(4753,-3382){\vector(-1,-2){855}}}{\put(5131,-3346){\line(1,0){270}}}{\put(5626,-3346){\line(1,0){270}}}{\put(6121,-3346){\line(1,0){270}}}{\put(6661,-3346){\line(1,0){270}}}{\put(7156,-3346){\line(1,0){270}}}{\put(7651,-3346){\line(1,0){270}}}{\put(4276,-5191){\line(1,0){270}}}{\put(4816,-5191){\line(1,0){270}}}{\put(5356,-5191){\line(1,0){270}}}{\put(5896,-5191){\line(1,0){270}}}{\put(6436,-5191){\line(1,0){270}}}{\put(6976,-5191){\line(1,0){270}}}{\put(7516,-5191){\line(1,0){270}}}{\put(2035,-5074){\vector(-1,2){855}}}{\put(2971,-3346){\vector(-1,-2){855}}}{\put(3943,-5110){\vector(-1,2){855}}}{\put(9073,-5110){\vector(-1,2){855}}}{\put(9946,-3391){\vector(-1,-2){855}}}{\put(10882,-5128){\vector(-1,2){855}}}{\put(10009,-3445){\vector(-1,-2){855}}}{\put(8551,-5191){\line(1,0){270}}}{\put(8056,-5191){\line(1,0){270}}}{\put(2161,-5191){\vector(1,0){1665}}}{\put(1261,-3346){\vector(1,0){1665}}}{\put(3016,-3346){\vector(1,0){1665}}}{\put(8281,-3346){\vector(1,0){1665}}}{\put(9181,-5191){\vector(1,0){1665}}}\put(746,-3166){\makebox(0,0)[lb]{\smash{{\SetFigFontNFSS{12}{14.4}{\rmdefault}{\mddefault}{\updefault}{$H$}}}}}\put(3221,-5506){\makebox(0,0)[lb]{\smash{{\SetFigFontNFSS{12}{14.4}{\rmdefault}{\mddefault}{\updefault}{$r_{n-2}r_{n-1}^2H$}}}}}\put(1421,-5506){\makebox(0,0)[lb]{\smash{{\SetFigFontNFSS{12}{14.4}{\rmdefault}{\mddefault}{\updefault}{$r_{n-1}^2H$}}}}}\put(2321,-3121){\makebox(0,0)[lb]{\smash{{\SetFigFontNFSS{12}{14.4}{\rmdefault}{\mddefault}{\updefault}{$r_{n-1}H$}}}}}\put(4031,-3121){\makebox(0,0)[lb]{\smash{{\SetFigFontNFSS{12}{14.4}{\rmdefault}{\mddefault}{\updefault}{$r_{n-3}r_{n-1}H$}}}}}\put(2456,-4651){\makebox(0,0)[lb]{\smash{{\SetFigFontNFSS{12}{14.4}{\rmdefault}{\mddefault}{\updefault}{$r_{n-2}$}}}}}\put(3356,-3976){\makebox(0,0)[lb]{\smash{{\SetFigFontNFSS{12}{14.4}{\rmdefault}{\mddefault}{\updefault}{$r_{n-3}$}}}}}\put(1511,-3976){\makebox(0,0)[lb]{\smash{{\SetFigFontNFSS{12}{14.4}{\rmdefault}{\mddefault}{\updefault}{$r_{n-1}$}}}}}\put(8666,-4021){\makebox(0,0)[lb]{\smash{{\SetFigFontNFSS{12}{14.4}{\rmdefault}{\mddefault}{\updefault}{$r_2$}}}}}\put(9566,-4741){\makebox(0,0)[lb]{\smash{{\SetFigFontNFSS{12}{14.4}{\rmdefault}{\mddefault}{\updefault}{$r_1$}}}}}\put(6821,-3121){\makebox(0,0)[lb]{\smash{{\SetFigFontNFSS{12}{14.4}{\rmdefault}{\mddefault}{\updefault}{$r_4\dots r_{n-3}r_{n-1}H$}}}}}\put(9296,-3121){\makebox(0,0)[lb]{\smash{{\SetFigFontNFSS{12}{14.4}{\rmdefault}{\mddefault}{\updefault}{$r_2r_4\dots r_{n-3}r_{n-1}H$}}}}}\put(9791,-5506){\makebox(0,0)[lb]{\smash{{\SetFigFontNFSS{12}{14.4}{\rmdefault}{\mddefault}{\updefault}{$r_1r_3\dots r_{n-2}r_{n-1}^2H$}}}}}\put(7136,-5506){\makebox(0,0)[lb]{\smash{{\SetFigFontNFSS{12}{14.4}{\rmdefault}{\mddefault}{\updefault}{$r_3\dots r_{n-2}r_{n-1}^2H$}}}}}
\put(1650,-6200){\begin{pic}\label{CoxToddaltA3b} Coxeter--Todd figure for $(A^+_n,H)$ for the presentation (\ref{codeAn}), $n$ odd
\end{pic}}
\end{picture}

\paragraph{4.} By induction on $n$, an upper bound, implied by the  Coxeter--Todd algorithm, 
for the cardinality of the group defined by (\ref{codeCar}) or (\ref{codeMoo})  or (\ref{codeAn}) is $n!/2$ in each case. The converse inequality 
is implied by the following surjective morphisms onto the group of even permutations of $n+1$ elements:
\begin{itemize}\item For Carmichael presentation (\ref{codeCar}): $a_i\mapsto(1,2,i+2)$, $i=1,\dots,n-1$. 
\item For Moore presentation (\ref{codeMoo}): $R_i\mapsto(1,2)(i+1,i+2)$, $i=1,\dots,n-1$. 
\item For the presentation (\ref{codeAn}): $r_i\mapsto(i,i+1,i+2)$, $i=1,\dots,n-1$. 
\end{itemize}
This concludes the proof that each of the groups -
defined by (\ref{codeCar}) or (\ref{codeMoo}) or (\ref{codeAn}) -
has the cardinality $n!/2$ and is isomorphic to the group $A^+_n$ (thus the morphisms above are the isomorphisms).

\bigskip
For each of the three presentations, the Coxeter--Todd figure gives a list $E_n$ of elements of $A^+_n$ such that any $x\in A^+_n$ can be written as $u_nh$ where $h\in A^+_{n-1}$ 
and $u_n\in E_n$; the set $E_n$ is in bijection with the set of vertices of the  Coxeter--Todd figure. 
The element $h$ in turn can be written in a normal form with respect to $A^+_{n-2}$.
Continuing we obtain recursively three normal forms for elements of $A^+_n$.

\begin{prop}
Let $A^+_n$ be given by Carmichael presentation 
(\ref{codeCar}). Any element $x\in A^+_n$ can be uniquely written as $x=u_nu_{n-1}\dots u_2$ where $u_i\in E_i$,  
\[E_i=\left\{1, a_{i-1},a_{i-1}^2,a_{i-2}a_{i-1}^2,\dots,a_2a_{i-1}^2,a_1a_{i-1}^2\right\}\ ,\ i=2,\dots,n.\]
\end{prop}

\begin{prop}
Let $A^+_n$ be given by Moore presentation (\ref{codeMoo}). Any element $x\in A^+_n$ can be uniquely written as 
$x=u_nu_{n-1}\dots u_2$ where $u_i\in E_i$, 
\[E_i=\left\{1,R_{i-1},R_{i-2}R_{i-1},\dots,R_1R_2\dots
R_{i-1},R_1^2R_2\dots R_{i-1}\right\}\ ,\ i=2,\dots,n.\]
\end{prop}

\begin{prop}\label{prop3} 
Let $A^+_n$ be given by the presentation (\ref{codeAn}). Any element $x\in A^+_n$ can be uniquely written as 
$x=u_nu_{n-1}\dots u_2$ where $u_i\in E_i$, $ i=2,\dots,n$, and  
\begin{itemize}
\item if $i$ is even, $E_i:=\left\{1,r_{i-1},r_{i-3}r_{i-1},\dots,r_1r_3\dots
r_{i-3}r_{i-1},r_{i-1}^2,r_{i-2}r_{i-1}^2,\dots,r_2r_4\dots r_{i-2}r_{i-1}^2\right\}$,
\item if $i$ is odd, $E_i:=\left\{1,r_{i-1},r_{i-3}r_{i-1},\dots,r_2r_4\dots
r_{i-3}r_{i-1},r_{i-1}^2,r_{i-2}r_{i-1}^2,\dots,r_1r_3\dots r_{i-2}r_{i-1}^2\right\}$.
\end{itemize}
\end{prop}

\subsection{Type B}

We label the vertices and oriented edges of the Coxeter graph of type B as follows.

\medskip
\setlength{\unitlength}{2000sp}
\begingroup\makeatletter\ifx\SetFigFontNFSS\undefined
\gdef\SetFigFontNFSS#1#2#3#4#5{
  \reset@font\fontsize{#1}{#2pt}
  \fontfamily{#3}\fontseries{#4}\fontshape{#5}
  \selectfont}
\fi\endgroup
\begin{picture}(11210,2000)(-800,-4779)
{\thinlines
\put(2989,-3346){\circle*{144}}}{\put(1171,-3346){\circle*{144}}}{\put(10081,-3346){\circle*{144}}}{\put(11926,-3346){\circle*{144}}}{\put(4861,-3346){\line(1,0){270}}}{\put(5356,-3346){\line(1,0){270}}}{\put(5851,-3346){\line(1,0){270}}}{\put(6346,-3346){\line(1,0){270}}}{\put(6886,-3346){\line(1,0){270}}}{\put(7381,-3346){\line(1,0){270}}}{\put(7876,-3346){\line(1,0){270}}}{\put(2881,-3346){\vector(-1,0){1665}}}{\put(4726,-3346){\vector(-1,0){1665}}}{\put(10036,-3346){\vector(-1,0){1665}}}{\put(11836,-3301){\vector(-1,0){1665}}}{\put(11836,-3391){\vector(-1,0){1665}}}\put(8756,-3256){\makebox(0,0)[lb]{\smash{{\SetFigFontNFSS{12}{14.4}{\rmdefault}{\mddefault}{\updefault}{$r_2$}}}}}\put(2546,-3706){\makebox(0,0)[lb]{\smash{{\SetFigFontNFSS{12}{14.4}{\rmdefault}{\mddefault}{\updefault}{$s_{n-2}$}}}}}\put(3446,-3256){\makebox(0,0)[lb]{\smash{{\SetFigFontNFSS{12}{14.4}{\rmdefault}{\mddefault}{\updefault}{$r_{n-2}$}}}}}\put(1556,-3256){\makebox(0,0)[lb]{\smash{{\SetFigFontNFSS{12}{14.4}{\rmdefault}{\mddefault}{\updefault}{$r_{n-1}$}}}}}\put(9656,-3706){\makebox(0,0)[lb]{\smash{{\SetFigFontNFSS{12}{14.4}{\rmdefault}{\mddefault}{\updefault}{$s_1$}}}}}\put(611,-3706){\makebox(0,0)[lb]{\smash{{\SetFigFontNFSS{12}{14.4}{\rmdefault}{\mddefault}{\updefault}{$s_{n-1}$}}}}}\put(11501,-3706){\makebox(0,0)[lb]{\smash{{\SetFigFontNFSS{12}{14.4}{\rmdefault}{\mddefault}{\updefault}{$s_0$}}}}}\put(10601,-3211){\makebox(0,0)[lb]{\smash{{\SetFigFontNFSS{12}{14.4}{\rmdefault}{\mddefault}{\updefault}{$r_1$}}}}}
\put(3400,-4300){\begin{pic}\label{Cox-B}Coxeter graph of type B
\end{pic}}
\end{picture}

\noindent The Coxeter group $B_n$ is generated by $s_0$, $\dots$, $s_{n-1}$ with the defining relations
\begin{equation}
\label{Bn}
\left\{\begin{array}{ll}
s_i^2=1 & \text{for $i=0,\dots,n-1$,}\\[.2em]
s_0s_1s_0s_1=s_1s_0s_1s_0, &\\[.2em]
s_is_{i+1}s_i=s_{i+1}s_is_{i+1} &\textrm{for $i=1,\dots,n-2$,}\\[.2em]
s_is_j=s_js_i & \textrm{for $i,j=0,\dots,n-1$ such that $|i-j|>1$.}
\end{array}\right.
\end{equation}

The group $B_n$ is isomorphic to the wreath product $C_2\wr A_{n-1}$ of the cyclic group $C_2$ of order $2$ by the symmetric group $A_{n-1}$. Denote by $\gamma$ the generator of $C_2$ and by $\eC$ the unit element of $C_2$. Let $\gamma^{(i)}:=(\eC,\dots,\eC,\gamma,\eC,\dots,\eC)$, $i=1,\dots,n$, be the element of $C_2^n$ with $\gamma$ in position $i$ and $\eC$ anywhere else; let $\eCn:=(\eC,\dots,\eC)$ denote the unit element of $C_2^n$ and $\eA$ the unit element of $A_{n-1}$. The isomorphism between $B_n$ and $C_2\wr A_{n-1}$ is given by $s_0\mapsto (\gamma^{(1)},\eA)$ and $s_i\mapsto \bigl(\eCn ,(i,i+1))$, $i=1,\dots,n-1$.

\vskip .2cm
Let $H$ be the subgroup generated by $s_0$, $\dots$, $s_{n-2}$; here is the  Coxeter--Todd figure for $(B_n,H)$:

\setlength{\unitlength}{2500sp}

\begingroup\makeatletter\ifx\SetFigFontNFSS\undefined
\gdef\SetFigFontNFSS#1#2#3#4#5{
  \reset@font\fontsize{#1}{#2pt}
  \fontfamily{#3}\fontseries{#4}\fontshape{#5}
  \selectfont}
\fi\endgroup
\begin{picture}(12110,1685)(600,-4270)
{\thinlines
\put(541,-3346){\circle*{144}}}{\put(1486,-3346){\circle*{144}}}{\put(2476,-3346){\circle*{144}}}{\put(4924,-3346){\circle*{144}}}{\put(5896,-3346){\circle*{144}}}{\put(6886,-3346){\circle*{144}}}{\put(7876,-3346){\circle*{144}}}{\put(10396,-3346){\circle*{144}}}{\put(11386,-3346){\circle*{144}}}{\put(12376,-3346){\circle*{144}}}{\put(586,-3346){\line(1,0){810}}}{\put(1576,-3346){\line(1,0){810}}}{\put(2566,-3346){\line(1,0){270}}}{\put(2971,-3346){\line(1,0){270}}}{\put(3376,-3346){\line(1,0){270}}}{\put(3781,-3346){\line(1,0){270}}}{\put(4186,-3346){\line(1,0){270}}}{\put(4591,-3346){\line(1,0){270}}}{\put(4996,-3346){\line(1,0){810}}}{\put(5986,-3346){\line(1,0){810}}}{\put(6976,-3346){\line(1,0){810}}}{\put(7966,-3346){\line(1,0){270}}}{\put(8371,-3346){\line(1,0){270}}}{\put(8776,-3346){\line(1,0){270}}}{\put(9226,-3346){\line(1,0){270}}}{\put(9676,-3346){\line(1,0){270}}}{\put(10081,-3346){\line(1,0){270}}}{\put(10486,-3346){\line(1,0){810}}}{\put(11476,-3346){\line(1,0){810}}}\put(111,-3706){\makebox(0,0)[lb]{\smash{{\SetFigFontNFSS{12}{14.4}{\rmdefault}{\mddefault}{\updefault}{$H$}}}}}\put(500,-3256){\makebox(0,0)[lb]{\smash{{\SetFigFontNFSS{12}{14.4}{\rmdefault}{\mddefault}{\updefault}{$s_{n-1}$}}}}}\put(1450,-3256){\makebox(0,0)[lb]{\smash{{\SetFigFontNFSS{12}{14.4}{\rmdefault}{\mddefault}{\updefault}{$s_{n-2}$}}}}}\put(4971,-3256){\makebox(0,0)[lb]{\smash{{\SetFigFontNFSS{12}{14.4}{\rmdefault}{\mddefault}{\updefault}{$s_1$}}}}}\put(5961,-3256){\makebox(0,0)[lb]{\smash{{\SetFigFontNFSS{12}{14.4}{\rmdefault}{\mddefault}{\updefault}{$s_0$}}}}}\put(6951,-3256){\makebox(0,0)[lb]{\smash{{\SetFigFontNFSS{12}{14.4}{\rmdefault}{\mddefault}{\updefault}{$s_1$}}}}}\put(10370,-3256){\makebox(0,0)[lb]{\smash{{\SetFigFontNFSS{12}{14.4}{\rmdefault}{\mddefault}{\updefault}{$s_{n-2}$}}}}}\put(11350,-3256){\makebox(0,0)[lb]{\smash{{\SetFigFontNFSS{12}{14.4}{\rmdefault}{\mddefault}{\updefault}{$s_{n-1}$}}}}}
\put(2000,-4100){\begin{pic}\label{CoxToddB} Coxeter--Todd figure for $(B_n,H)$ for the presentation (\ref{Bn})
\end{pic}}
\end{picture}

\vskip .2cm
This is an analogue - for the group $B_n^+$ - of the Carmichael presentation (the distinguished oriented edge is $(0,1)$): 
\begin{itemize}
\item the alternating group $B^+_n$ is generated by $a_1$, $\dots$, $a_{n-1}$ with the defining relations
\begin{equation}\label{codeBn3}\left\{\begin{array}{ll}
a_i^4=1 & \text{for $i=1,\dots,n-1$,}\\[0.2em]
(a_1a_i)^3=1 & \text{for $i=2,\dots,n-1$,}\\[0.2em]
(a_1^2a_i)^2=1 & \text{for $i=2,\dots,n-1$.}\\[0.2em]
(a_1a_ia_1a_j)^2=1 & \text{for $i,j=2,\dots,n-1$ such that $i<j$.}
\end{array}\right.
\end{equation}
\end{itemize}
The presentation (\ref{codebourb}), respectively the presentation of the Proposition \ref{prop-code}, reads:
\begin{itemize}
 \item the alternating group $B^+_n$ is generated by $R_1$, $\dots$, $R_{n-1}$ with the defining relations
\begin{equation}\label{codeBn1}\left\{\begin{array}{ll}
R_1^4=1,\\[.2em]
R_i^2=1 & \textrm{for $i=2,\dots,n-1$,}\\[.2em]
(R_i^{-1}R_{i+1})^3=1 & \text{for $i=1,\dots,n-2$,}\\[.2em]
(R_i^{-1}R_j)^2=1 & \textrm{for $i,j=1,\dots,n-1$ such that $|i-j|>1$;}
\end{array}\right.
\end{equation}
\item respectively, the alternating group $B^+_n$ is generated by $r_1$, $\dots$, $r_{n-1}$ with the defining relations
\begin{equation}\label{codeBn2}\left\{\begin{array}{ll}
r_1^4=1,\\[.2em]
r_i^3=1 & \textrm{for $i=2,\dots,n-1$,}\\[.2em]
(r_ir_{i+1})^2=1 & \text{for $i=1,\dots,n-2$,}\\[.2em]
(r_ir_{i+1}r_{i+2})^2=1 & \text{for $i=1,\dots,n-3$,}\\[.2em]
r_ir_j=r_jr_i & \textrm{for $i,j=1,\dots,n-1$ such that $|i-j|>2$.}
\end{array}\right.
\end{equation}
\end{itemize}

Denote by $\epsilon_0$ the sign character of $A_{n-1}$. Let $\epsCn\ \colon\ C_2^n\to\{-1,1\}$ be the homomorphism defined by $\epsCn(\gamma^{(i)})=-1$ for $i=1,\dots,n$. The alternating group $B^+_n$ of type B is isomorphic to the subgroup of $C_2\wr A_{n-1}$ formed by elements $(g,\pi)$, $\pi\in A_{n-1}$ and $g\in C_2^n$, such that $\epsCn(g)\epsilon_0(\pi)=1$.
The isomorphisms are given by:
\begin{itemize}
\item for the presentation (\ref{codeBn3}), $a_i\mapsto\bigl(\gamma^{(1)},(1,i+1)\bigr)$, $i=1,\dots,n-1$.
\item for the presentation (\ref{codeBn1}), $R_i\mapsto\bigl(\gamma^{(1)},(i,i+1)\bigr)$, $i=1,\dots,n-1$.
\item for the presentation (\ref{codeBn2}), $r_1\mapsto\bigl(\gamma^{(1)},(1,2)\bigr)$ and $r_i\mapsto\bigl(\eCn ,(i-1,i,i+1)\bigr)$, $i=2,\dots,n-1$.
\end{itemize}

Let $H$ be the subgroup generated by $a_1,\dots,a_{n-2}$, or by $R_1$, $\dots$, $R_{n-2}$, or by $r_1$, $\dots$, $r_{n-2}$.

\vskip 1cm
\setlength{\unitlength}{2450sp}
\begingroup\makeatletter\ifx\SetFigFontNFSS\undefined
\gdef\SetFigFontNFSS#1#2#3#4#5{
  \reset@font\fontsize{#1}{#2pt}
  \fontfamily{#3}\fontseries{#4}\fontshape{#5}
  \selectfont}
\fi\endgroup
\begin{picture}(2283,4528)(371,-6560)
{\put(5300,-3100){\circle*{120}}}
\put(170,-2850){\makebox(0,0)[lb]{\smash{{\SetFigFontNFSS{12}{14.4}{\rmdefault}{\mddefault}{\updefault}{$H$}}}}}
{\put(5300,-3100){\vector(2,3){900}}}
{\put(6260,-1720){\circle*{120}}}
{\put(6320,-1750){\vector(2,-3){880}}}
{\put(7220,-3100){\circle*{120}}}
{\put(7200,-3150){\vector(-2,-3){870}}}
{\put(6300,-4480){\circle*{120}}}
{\put(6250,-4480){\vector(-2,3){900}}}
\put(5850,-3150){\makebox(0,0)[lb]{\smash{{\SetFigFontNFSS{12}{14.4}{\rmdefault}{\mddefault}{\updefault}{$a_1$}}}}}
{\put(3450,-3100){\circle*{120}}}
{\put(3500,-3050){\vector(2,1){2730}}}
{\put(6320,-1720){\vector(2,-1){2730}}}
{\put(9100,-3100){\circle*{120}}}
{\put(9030,-3140){\vector(-2,-1){2710}}}
{\put(6230,-4480){\vector(-2,1){2750}}}
\put(3730,-3150){\makebox(0,0)[lb]{\smash{{\SetFigFontNFSS{12}{14.4}{\rmdefault}{\mddefault}{\updefault}{$a_2$}}}}}
\put(8000,-3150){\makebox(0,0)[lb]{\smash{{\SetFigFontNFSS{12}{14.4}{\rmdefault}{\mddefault}{\updefault}{$a_2$}}}}}
{\put(550,-3100){\circle*{120}}}
{\put(550,-3080){\vector(4,1){5640}}}
{\put(6300,-1680){\vector(4,-1){5620}}}
{\put(11980,-3100){\circle*{120}}}
{\put(11940,-3110){\vector(-4,-1){5580}}}
{\put(6240,-4500){\vector(-4,1){5630}}}
\put(850,-3150){\makebox(0,0)[lb]{\smash{{\SetFigFontNFSS{12}{14.4}{\rmdefault}{\mddefault}{\updefault}{$a_{n-1}\ \ \dots\ \ \dots$}}}}}
\put(9300,-3150){\makebox(0,0)[lb]{\smash{{\SetFigFontNFSS{12}{14.4}{\rmdefault}{\mddefault}{\updefault}{$\dots\ \ \dots\ \ a_{n-1}$}}}}}
\put(2000,-5250){\begin{pic}\label{CoxToddaltB3} Coxeter--Todd figure for $(B^+_n,H)$ for the presentation (\ref{codeBn3})
\end{pic}}
\end{picture}

\vspace{-1.2cm}
\setlength{\unitlength}{2500sp}
\begingroup\makeatletter\ifx\SetFigFontNFSS\undefined
\gdef\SetFigFontNFSS#1#2#3#4#5{
  \reset@font\fontsize{#1}{#2pt}
  \fontfamily{#3}\fontseries{#4}\fontshape{#5}
  \selectfont}
\fi\endgroup
\begin{picture}(12110,2600)(346,-4595)
{\thinlines
\put(541,-3346){\circle*{144}}}{\put(1486,-3346){\circle*{144}}}{\put(2476,-3346){\circle*{144}}}{\put(6886,-3346){\circle*{144}}}{\put(7876,-3346){\circle*{144}}}{\put(10396,-3346){\circle*{144}}}{\put(11386,-3346){\circle*{144}}}{\put(12376,-3346){\circle*{144}}}{\put(3736,-3346){\circle*{144}}}{\put(4726,-3346){\circle*{144}}}{\put(5806,-2266){\circle*{144}}}{\put(5824,-4516){\circle*{144}}}{\put(586,-3346){\line(1,0){810}}}{\put(1576,-3346){\line(1,0){810}}}{\put(2566,-3346){\line(1,0){270}}}{\put(2971,-3346){\line(1,0){270}}}{\put(3376,-3346){\line(1,0){270}}}{\put(6976,-3346){\line(1,0){810}}}{\put(7966,-3346){\line(1,0){270}}}{\put(8371,-3346){\line(1,0){270}}}{\put(8776,-3346){\line(1,0){270}}}{\put(9226,-3346){\line(1,0){270}}}{\put(9676,-3346){\line(1,0){270}}}{\put(10081,-3346){\line(1,0){270}}}{\put(10486,-3346){\line(1,0){810}}}{\put(11476,-3346){\line(1,0){810}}}{\put(3826,-3346){\line(1,0){810}}}{\put(4748,-3278){\vector(1,1){1035}}}{\put(5851,-2266){\vector(1,-1){1035}}}{\put(6930,-3391){\vector(-1,-1){1035}}}{\put(5761,-4471){\vector(-1,1){1035}}}{\put(5896,-2311){\line(0,-1){2160}}}{\put(5806,-4471){\line(0, 1){2115}}}{\put(5716,-4471){\line(0, 1){2115}}}\put(111,-3706){\makebox(0,0)[lb]{\smash{{\SetFigFontNFSS{12}{14.4}{\rmdefault}{\mddefault}{\updefault}{\small{$H$}}}}}}\put(6951,-3256){\makebox(0,0)[lb]{\smash{{\SetFigFontNFSS{12}{14.4}{\rmdefault}{\mddefault}{\updefault}{\small{$R_2$}}}}}}\put(3711,-3256){\makebox(0,0)[lb]{\smash{{\SetFigFontNFSS{12}{14.4}{\rmdefault}{\mddefault}{\updefault}{\small{$R_2$}}}}}}\put(4666,-2716){\makebox(0,0)[lb]{\smash{{\SetFigFontNFSS{12}{14.4}{\rmdefault}{\mddefault}{\updefault}{\small{$R_1$}}}}}}\put(381,-3256){\makebox(0,0)[lb]{\smash{{\SetFigFontNFSS{12}{14.4}{\rmdefault}{\mddefault}{\updefault}{\small{$R_{n-1}$}}}}}}\put(1326,-3256){\makebox(0,0)[lb]{\smash{{\SetFigFontNFSS{12}{14.4}{\rmdefault}{\mddefault}{\updefault}{\small{$R_{n-2}$}}}}}}\put(10281,-3256){\makebox(0,0)[lb]{\smash{{\SetFigFontNFSS{12}{14.4}{\rmdefault}{\mddefault}{\updefault}{\small{$R_{n-2}$}}}}}}\put(11300,-3256){\makebox(0,0)[lb]{\smash{{\SetFigFontNFSS{12}{14.4}{\rmdefault}{\mddefault}{\updefault}{\small{$R_{n-1}$}}}}}}\put(4890,-3436){\makebox(0,0)[lb]{\smash{{\SetFigFontNFSS{12}{14.4}{\rmdefault}{\mddefault}{\updefault}{\small{$R_2,\dots,R_{n-1}$}}}}}}\put(4666,-4201){\makebox(0,0)[lb]{\smash{{\SetFigFontNFSS{12}{14.4}{\rmdefault}{\mddefault}{\updefault}{\small{$R_1$}}}}}}\put(6141,-4201){\makebox(0,0)[lb]{\smash{{\SetFigFontNFSS{12}{14.4}{\rmdefault}{\mddefault}{\updefault}{\small{$R_1$}}}}}}\put(6051,-2716){\makebox(0,0)[lb]{\smash{{\SetFigFontNFSS{12}{14.4}{\rmdefault}{\mddefault}{\updefault}{\small{$R_1$}}}}}}
\put(2000,-5000){\begin{pic}\label{CoxToddaltB1} Coxeter--Todd figure for $(B^+_n,H)$ for the presentation (\ref{codeBn1})
\end{pic}}
\end{picture}

\bigskip
\setlength{\unitlength}{2400sp}

\begingroup\makeatletter\ifx\SetFigFontNFSS\undefined
\gdef\SetFigFontNFSS#1#2#3#4#5{
  \reset@font\fontsize{#1}{#2pt}
  \fontfamily{#3}\fontseries{#4}\fontshape{#5}
  \selectfont}
\fi\endgroup
\begin{picture}(12515,2700)(76,-5225)
{\thinlines
\put(11611,-3346){\circle*{144}}}{\put(4996,-3346){\circle*{144}}}{\put(1171,-5146){\circle*{144}}}{\put(2926,-5146){\circle*{144}}}{\put(271,-3346){\circle*{144}}}{\put(2098,-3346){\circle*{144}}}{\put(4168,-5146){\circle*{144}}}{\put(5896,-5146){\circle*{144}}}{\put(6841,-3346){\circle*{144}}}{\put(7696,-5146){\circle*{144}}}{\put(8582,-3365){\circle*{144}}}{\put(9568,-5146){\circle*{144}}}{\put(10711,-5146){\circle*{144}}}{\put(12511,-5101){\circle*{144}}}{\put(4978,-3427){\vector(-1,-2){855}}}{\put(1126,-5146){\vector(-1,2){855}}}{\put(2116,-3391){\vector(-1,-2){855}}}{\put(2053,-3337){\vector(-1,-2){855}}}{\put(2980,-5119){\vector(-1,2){855}}}{\put(2341,-3346){\line(1,0){270}}}{\put(2791,-3346){\line(1,0){270}}}{\put(3241,-3346){\line(1,0){270}}}{\put(3691,-3346){\line(1,0){270}}}{\put(4591,-3346){\line(1,0){270}}}{\put(4141,-3346){\line(1,0){270}}}{\put(3691,-5146){\line(1,0){270}}}{\put(3196,-5146){\line(1,0){270}}}{\put(5824,-5137){\vector(-1,2){855}}}{\put(5914,-5092){\vector(-1,2){855}}}{\put(6796,-3436){\vector(1,-2){855}}}{\put(6868,-3400){\vector(1,-2){855}}}{\put(7678,-5092){\vector(1,2){855}}}{\put(7759,-5155){\vector(1,2){855}}}{\put(9523,-5110){\vector(-1,2){855}}}{\put(8776,-3346){\line(1,0){270}}}{\put(9226,-3346){\line(1,0){270}}}{\put(9721,-5146){\line(1,0){270}}}{\put(10171,-5146){\line(1,0){270}}}{\put(11629,-3337){\vector(1,-2){855}}}{\put(10729,-5065){\vector(1,2){855}}}{\put(9721,-3346){\line(1,0){270}}}{\put(10171,-3346){\line(1,0){270}}}{\put(10621,-3346){\line(1,0){270}}}{\put(11116,-3346){\line(1,0){270}}}{\put(1261,-5146){\vector(1,0){1665}}}{\put(361,-3346){\vector(1,0){1665}}}{\put(4186,-5146){\vector(1,0){1665}}}{\put(7606,-5146){\vector(-1,0){1665}}}{\put(5086,-3346){\vector(1,0){1665}}}{\put(8596,-3346){\vector(-1,0){1665}}}{\put(9451,-5146){\vector(-1,0){1665}}}{\put(12466,-5146){\vector(-1,0){1665}}}\put(521,-3931){\makebox(0,0)[lb]{\smash{{\SetFigFontNFSS{12}{14.4}{\rmdefault}{\mddefault}{\updefault}{$r_{n-1}$}}}}}\put(-109,-3166){\makebox(0,0)[lb]{\smash{{\SetFigFontNFSS{12}{14.4}{\rmdefault}{\mddefault}{\updefault}{$H$}}}}}\put(1601,-4696){\makebox(0,0)[lb]{\smash{{\SetFigFontNFSS{12}{14.4}{\rmdefault}{\mddefault}{\updefault}{$r_{n-2}$}}}}}\put(4526,-4696){\makebox(0,0)[lb]{\smash{{\SetFigFontNFSS{12}{14.4}{\rmdefault}{\mddefault}{\updefault}{$r_2$}}}}}\put(5786,-4291){\makebox(0,0)[lb]{\smash{{\SetFigFontNFSS{12}{14.4}{\rmdefault}{\mddefault}{\updefault}{$r_1$}}}}}\put(7226,-3976){\makebox(0,0)[lb]{\smash{{\SetFigFontNFSS{12}{14.4}{\rmdefault}{\mddefault}{\updefault}{$r_2$}}}}}\put(8216,-4741){\makebox(0,0)[lb]{\smash{{\SetFigFontNFSS{12}{14.4}{\rmdefault}{\mddefault}{\updefault}{$r_3$}}}}}\put(11096,-4741){\makebox(0,0)[lb]{\smash{{\SetFigFontNFSS{12}{14.4}{\rmdefault}{\mddefault}{\updefault}{$r_{n-1}$}}}}}
\put(1200,-5800){\begin{pic}\label{CoxToddaltB2a} Coxeter--Todd figure for $(B^+_n,H)$ for the presentation (\ref{codeBn2}), $n$ even
\end{pic}}
\end{picture}

\bigskip
\setlength{\unitlength}{2400sp}

\begingroup\makeatletter\ifx\SetFigFontNFSS\undefined
\gdef\SetFigFontNFSS#1#2#3#4#5{
  \reset@font\fontsize{#1}{#2pt}
  \fontfamily{#3}\fontseries{#4}\fontshape{#5}
  \selectfont}
\fi\endgroup
\begin{picture}(12020,3300)(571,-5770)
{\thinlines
\put(811,-3346){\circle*{144}}}{\put(1711,-5146){\circle*{144}}}{\put(2611,-3346){\circle*{144}}}{\put(3511,-5146){\circle*{144}}}{\put(4591,-3346){\circle*{144}}}{\put(5491,-5191){\circle*{144}}}{\put(6391,-3346){\circle*{144}}}{\put(7246,-5191){\circle*{144}}}{\put(8119,-3346){\circle*{144}}}{\put(9046,-5191){\circle*{144}}}{\put(9901,-3346){\circle*{144}}}{\put(12511,-5191){\circle*{144}}}{\put(10756,-5191){\circle*{144}}}{\put(11611,-3346){\circle*{144}}}{\put(1675,-5119){\vector(-1,2){855}}}{\put(2548,-3382){\vector(-1,-2){855}}}{\put(2656,-3391){\vector(-1,-2){855}}}{\put(3520,-5119){\vector(-1,2){855}}}{\put(2746,-3346){\line(1,0){270}}}{\put(3196,-3346){\line(1,0){270}}}{\put(3646,-3346){\line(1,0){270}}}{\put(4096,-3346){\line(1,0){270}}}{\put(3646,-5146){\line(1,0){270}}}{\put(5482,-5128){\vector(-1,2){855}}}{\put(4951,-5146){\line(1,0){270}}}{\put(4051,-5146){\line(1,0){270}}}{\put(4501,-5146){\line(1,0){270}}}{\put(6310,-3418){\vector(-1,-2){855}}}{\put(6400,-3463){\vector(-1,-2){855}}}{\put(9055,-5128){\vector(1,2){855}}}{\put(9991,-3346){\line(1,0){270}}}{\put(10396,-3346){\line(1,0){270}}}{\put(10801,-3346){\line(1,0){270}}}{\put(11206,-3346){\line(1,0){270}}}{\put(11629,-3427){\vector(1,-2){855}}}{\put(9181,-5191){\line(1,0){270}}}{\put(9586,-5191){\line(1,0){270}}}{\put(9946,-5191){\line(1,0){270}}}{\put(10351,-5191){\line(1,0){270}}}{\put(10729,-5110){\vector(1,2){855}}}{\put(8200,-3409){\vector(1,-2){855}}}{\put(8110,-3454){\vector(1,-2){855}}}{\put(7318,-5092){\vector(1,2){855}}}{\put(7201,-5101){\vector(1,2){855}}}{\put(901,-3346){\vector(1,0){1665}}}{\put(1801,-5146){\vector(1,0){1665}}}{\put(4681,-3346){\vector(1,0){1665}}}{\put(5581,-5191){\vector(1,0){1665}}}{\put(8101,-3346){\vector(-1,0){1665}}}{\put(9001,-5191){\vector(-1,0){1665}}}{\put(9811,-3346){\vector(-1,0){1665}}}{\put(12466,-5191){\vector(-1,0){1665}}}\put(386,-3166){\makebox(0,0)[lb]{\smash{{\SetFigFontNFSS{12}{14.4}{\rmdefault}{\mddefault}{\updefault}{$H$}}}}}\put(1106,-3976){\makebox(0,0)[lb]{\smash{{\SetFigFontNFSS{12}{14.4}{\rmdefault}{\mddefault}{\updefault}{$r_{n-1}$}}}}}\put(5111,-4021){\makebox(0,0)[lb]{\smash{{\SetFigFontNFSS{12}{14.4}{\rmdefault}{\mddefault}{\updefault}{$r_2$}}}}}\put(2096,-4651){\makebox(0,0)[lb]{\smash{{\SetFigFontNFSS{12}{14.4}{\rmdefault}{\mddefault}{\updefault}{$r_{n-2}$}}}}}\put(6371,-4336){\makebox(0,0)[lb]{\smash{{\SetFigFontNFSS{12}{14.4}{\rmdefault}{\mddefault}{\updefault}{$r_1$}}}}}\put(8576,-4021){\makebox(0,0)[lb]{\smash{{\SetFigFontNFSS{12}{14.4}{\rmdefault}{\mddefault}{\updefault}{$r_3$}}}}}\put(7676,-4741){\makebox(0,0)[lb]{\smash{{\SetFigFontNFSS{12}{14.4}{\rmdefault}{\mddefault}{\updefault}{$r_2$}}}}}\put(11141,-4696){\makebox(0,0)[lb]{\smash{{\SetFigFontNFSS{12}{14.4}{\rmdefault}{\mddefault}{\updefault}{$r_{n-1}$}}}}}
\put(1700,-5800){\begin{pic}\label{CoxToddaltB2b} Coxeter--Todd figure for $(B^+_n,H)$ for the presentation (\ref{codeBn2}), $n$ odd
\end{pic}}
\end{picture}

\bigskip
As for type A, the  Coxeter--Todd figures for the presentations (\ref{codeBn3}), (\ref{codeBn1}) and (\ref{codeBn2}) provide three normal forms for elements of the alternating group $B^+_n$. We omit details.

\vskip .2cm
\textbf{Remark.} The complex reflection group $G(m,1,n)$ is generated by $s_0$, $\dots$, $s_{n-1}$ with the same defining relations as in (\ref{Bn}), except that $s_0^2$ is replaced by $s_0^m=1$. For $m=1$,  $G(1,1,n)$ is the group $A_{n-1}$; for $m=2$, $G(2,1,n)$ is the group $B_n$. Let $H$ be the subgroup of $G(m,1,n)$ generated by $s_0$, $\dots$, $s_{n-2}$.
The Coxeter-Todd figure for $\bigl(G(m,1,n), H\bigr)$, see \cite{OPdA2}, generalizes Fig. \ref{CoxToddA} and Fig. \ref{CoxToddB}.

\subsection{Type D}

We label the vertices and oriented edges of the Coxeter graph of type D as follows.

\bigskip
\setlength{\unitlength}{1800sp}
\begingroup\makeatletter\ifx\SetFigFontNFSS\undefined
\gdef\SetFigFontNFSS#1#2#3#4#5{
  \reset@font\fontsize{#1}{#2pt}
  \fontfamily{#3}\fontseries{#4}\fontshape{#5}
  \selectfont}
\fi\endgroup
\begin{picture}(10740,3228)(-1400,-5155)
{\thinlines
\put(2989,-3346){\circle*{144}}}{\put(1171,-3346){\circle*{144}}}{\put(10081,-3346){\circle*{144}}}{\put(11386,-2086){\circle*{144}}}{\put(11386,-4651){\circle*{144}}}{\put(4861,-3346){\line(1,0){270}}}{\put(5356,-3346){\line(1,0){270}}}{\put(5851,-3346){\line(1,0){270}}}{\put(6346,-3346){\line(1,0){270}}}{\put(6886,-3346){\line(1,0){270}}}{\put(7381,-3346){\line(1,0){270}}}{\put(7876,-3346){\line(1,0){270}}}{\put(2881,-3346){\vector(-1,0){1665}}}{\put(4726,-3346){\vector(-1,0){1665}}}{\put(10036,-3346){\vector(-1,0){1665}}}{\put(11324,-2102){\vector(-1,-1){1177}}}{\put(11320,-4585){\vector(-1,1){1177}}}\put(8756,-3256){\makebox(0,0)[lb]{\smash{{\SetFigFontNFSS{12}{14.4}{\rmdefault}{\mddefault}{\updefault}{$r_3$}}}}}\put(2346,-3706){\makebox(0,0)[lb]{\smash{{\SetFigFontNFSS{12}{14.4}{\rmdefault}{\mddefault}{\updefault}{$s_{n-2}$}}}}}\put(3246,-3256){\makebox(0,0)[lb]{\smash{{\SetFigFontNFSS{12}{14.4}{\rmdefault}{\mddefault}{\updefault}{$r_{n-2}$}}}}}\put(1356,-3256){\makebox(0,0)[lb]{\smash{{\SetFigFontNFSS{12}{14.4}{\rmdefault}{\mddefault}{\updefault}{$r_{n-1}$}}}}}\put(9456,-3706){\makebox(0,0)[lb]{\smash{{\SetFigFontNFSS{12}{14.4}{\rmdefault}{\mddefault}{\updefault}{$s_2$}}}}}\put(411,-3706){\makebox(0,0)[lb]{\smash{{\SetFigFontNFSS{12}{14.4}{\rmdefault}{\mddefault}{\updefault}{$s_{n-1}$}}}}}\put(10106,-2536){\makebox(0,0)[lb]{\smash{{\SetFigFontNFSS{12}{14.4}{\rmdefault}{\mddefault}{\updefault}{$r_2$}}}}}\put(10151,-4291){\makebox(0,0)[lb]{\smash{{\SetFigFontNFSS{12}{14.4}{\rmdefault}{\mddefault}{\updefault}{$r_1$}}}}}\put(11121,-2086){\makebox(0,0)[lb]{\smash{{\SetFigFontNFSS{12}{14.4}{\rmdefault}{\mddefault}{\updefault}{$s_0$}}}}}\put(11121,-4786){\makebox(0,0)[lb]{\smash{{\SetFigFontNFSS{12}{14.4}{\rmdefault}{\mddefault}{\updefault}{$s_1$}}}}}
\put(2500,-5250){\begin{pic}\label{Cox-D}Coxeter graph of type D
\end{pic}}
\end{picture}

\bigskip
The Coxeter group $D_n$ is generated by $s_0$, $\dots$, $s_{n-1}$ with the defining relations
\begin{equation}
\label{Dn}
\left\{\begin{array}{l}
s_i^2=1\ \ \ \text{for $i=0,\dots,n-1$,}\\[.2em]
s_0s_2s_0=s_2s_0s_2,\ \ s_is_{i+1}s_i=s_{i+1}s_is_{i+1}\ \ \ \textrm{for $i=1,\dots,n-2$,}\\[.2em]
s_0s_1=s_1s_0,\ \ s_0s_i=s_is_0\ \ \ \textrm{for $i=3,\dots,n-1$,}\\[.2em]
s_is_j=s_js_i\ \ \ \textrm{for $i,j=1,\dots,n-1$ such that $|i-j|>1$.}
\end{array}\right.
\end{equation}

The group $D_n$ is isomorphic to the subgroup of $C_2\wr A_{n-1}$ formed by elements $(g,\pi )$, $\pi\in A_{n-1}$ and $g\in C_2^n$, such that $\epsCn(g)=1$. The isomorphism is given by $s_0\mapsto\bigl(\gamma^{(1)}\gamma^{(2)},(1,2)\bigr)$ and $s_i\mapsto\bigl(\eCn ,(i,i+1)\bigr)$ for $i=1,\dots,n-1$.

\vskip .2cm
Let $H$ be the subgroup generated by $s_0$, $\dots$, $s_{n-2}$; here is the Coxeter--Todd figure for $(D_n,H)$:

\bigskip
\setlength{\unitlength}{2500sp}

\begingroup\makeatletter\ifx\SetFigFontNFSS\undefined
\gdef\SetFigFontNFSS#1#2#3#4#5{
  \reset@font\fontsize{#1}{#2pt}
  \fontfamily{#3}\fontseries{#4}\fontshape{#5}
  \selectfont}
\fi\endgroup
\begin{picture}(12110,2783)(346,-4790)
{\thinlines
\put(541,-3346){\circle*{144}}}{\put(1486,-3346){\circle*{144}}}{\put(2476,-3346){\circle*{144}}}{\put(10396,-3346){\circle*{144}}}{\put(11386,-3346){\circle*{144}}}{\put(12376,-3346){\circle*{144}}}{\put(4186,-3346){\circle*{144}}}{\put(5176,-3346){\circle*{144}}}{\put(6481,-2086){\circle*{144}}}{\put(8686,-3346){\circle*{144}}}{\put(7696,-3346){\circle*{144}}}{\put(6436,-4561){\circle*{144}}}{\put(586,-3346){\line(1,0){810}}}{\put(1576,-3346){\line(1,0){810}}}{\put(2566,-3346){\line(1,0){270}}}{\put(2971,-3346){\line(1,0){270}}}{\put(3376,-3346){\line(1,0){270}}}{\put(3781,-3346){\line(1,0){270}}}{\put(8776,-3346){\line(1,0){270}}}{\put(9226,-3346){\line(1,0){270}}}{\put(9676,-3346){\line(1,0){270}}}{\put(10081,-3346){\line(1,0){270}}}{\put(10486,-3346){\line(1,0){810}}}{\put(11476,-3346){\line(1,0){810}}}{\put(4276,-3346){\line(1,0){810}}}{\put(7786,-3346){\line(1,0){810}}}{\put(5176,-3346){\line(1,1){1260}}}{\put(6526,-2131){\line(1,-1){1125}}}{\put(5243,-3413){\line(1,-1){1125}}}{\put(6458,-4583){\line(1,1){1260}}}\put(61,-3706){\makebox(0,0)[lb]{\smash{{\SetFigFontNFSS{12}{14.4}{\rmdefault}{\mddefault}{\updefault}{$H$}}}}}\put(456,-3256){\makebox(0,0)[lb]{\smash{{\SetFigFontNFSS{12}{14.4}{\rmdefault}{\mddefault}{\updefault}{$s_{n-1}$}}}}}\put(1356,-3256){\makebox(0,0)[lb]{\smash{{\SetFigFontNFSS{12}{14.4}{\rmdefault}{\mddefault}{\updefault}{$s_{n-2}$}}}}}\put(10311,-3256){\makebox(0,0)[lb]{\smash{{\SetFigFontNFSS{12}{14.4}{\rmdefault}{\mddefault}{\updefault}{$s_{n-2}$}}}}}\put(11301,-3256){\makebox(0,0)[lb]{\smash{{\SetFigFontNFSS{12}{14.4}{\rmdefault}{\mddefault}{\updefault}{$s_{n-1}$}}}}}\put(7750,-3256){\makebox(0,0)[lb]{\smash{{\SetFigFontNFSS{12}{14.4}{\rmdefault}{\mddefault}{\updefault}{$s_2$}}}}}\put(4256,-3256){\makebox(0,0)[lb]{\smash{{\SetFigFontNFSS{12}{14.4}{\rmdefault}{\mddefault}{\updefault}{$s_2$}}}}}\put(5350,-4291){\makebox(0,0)[lb]{\smash{{\SetFigFontNFSS{12}{14.4}{\rmdefault}{\mddefault}{\updefault}{$s_1$}}}}}\put(6766,-4246){\makebox(0,0)[lb]{\smash{{\SetFigFontNFSS{12}{14.4}{\rmdefault}{\mddefault}{\updefault}{$s_0$}}}}}\put(5301,-2581){\makebox(0,0)[lb]{\smash{{\SetFigFontNFSS{12}{14.4}{\rmdefault}{\mddefault}{\updefault}{$s_0$}}}}}\put(6721,-2536){\makebox(0,0)[lb]{\smash{{\SetFigFontNFSS{12}{14.4}{\rmdefault}{\mddefault}{\updefault}{$s_1$}}}}}
\put(2100,-5000){\begin{pic}\label{CoxToddD} Coxeter--Todd figure for $(D_n,H)$ for the presentation (\ref{Dn})
\end{pic}}
\end{picture}

\vspace{.6cm}
This is an analogue - for the group $D_n^+$ with $n\geq 3$ - of the Carmichael presentation (the distinguished oriented edge is $(0,2)$): 
\begin{itemize}
\item the alternating group $D^+_n$, $n\geq 3$, is generated by $a_1$, $\dots$, $a_{n-1}$ with the defining relations
\begin{equation}\label{codeDn3}\left\{\begin{array}{ll}
a_i^3=1 & \text{for $i=1,\dots,n-1$,}\\[0.2em]
(a_1a_i)^2=1 & \text{for $i=2,\dots,n-1$,}\\[0.2em]
(a_2^2a_i)^2=1 & \text{for $i=3,\dots,n-1$,}\\[0.2em]
(a_ia_j)^2=1 & \text{for $i,j=3,\dots,n-1$ such that $i<j$.}
\end{array}\right.
\end{equation}
\end{itemize}
The presentation (\ref{codebourb}), respectively the presentation of the Proposition \ref{prop-code},
reads:
\begin{itemize}
 \item the alternating group $D^+_n$ is generated by $R_1$, $\dots$, $R_{n-1}$ with the defining relations
\begin{equation}\label{codeDn1}\left\{\begin{array}{ll}
R_2^3=1,\\[.2em]
R_i^2=1 & \textrm{for $i=1,\dots,n-1$ and $i\neq2$,}\\[.2em]
(R_i^{-1}R_{i+1})^3=1 & \text{for $i=1,\dots,n-2$,}\\[.2em]
(R_i^{-1}R_j)^2=1 & \textrm{for $i,j=1,\dots,n-1$ such that $|i-j|>1$;}
\end{array}\right.
\end{equation}
\item respectively, the alternating group $D^+_n$, $n\geq 3$, is generated by $r_1$, $\dots$, $r_{n-1}$ with the defining relations
\begin{equation}\label{codeDn2}\left\{\begin{array}{l}
r_i^3=1\ \ \ \text{for $i=1,\dots,n-1$,}\\[.2em]
(r_1r_2^2)^2=1,\ \ (r_1r_3)^2=1,\ \ (r_ir_{i+1})^2=1\ \ \ \textrm{for $i=2,\dots,n-2$,}\\[.2em]
(r_1r_3r_4)^2=1,\ \ (r_ir_{i+1}r_{i+2})^2=1\ \ \ \textrm{for $i=2,\dots,n-3$,}\\[.2em]
r_1r_i=r_ir_1\ \ \ \textrm{for $i=5,\dots,n-1$,}\\[.2em]
r_ir_j=r_jr_i\ \ \ \textrm{for $i,j=2,\dots,n-1$ such that $|i-j|>2$.}
\end{array}\right.
\end{equation}
\end{itemize}
The alternating group $D^+_n$ of type D is isomorphic to the subgroup of $C_2\wr A_{n-1}$ formed by elements $(g,\pi )$, $\pi\in A_{n-1}$ and $g\in C_2^n$, such that $\epsCn(g)=1$ and $\epsilon_0(\pi)=1$. The isomorphisms are given by:
\begin{itemize}
\item for the presentation (\ref{codeDn3}), 

$a_1\mapsto\bigl(\gamma^{(1)}\gamma^{(2)},(1,2,3)\bigr)$, $a_2\mapsto\bigl(\gamma^{(1)}\gamma^{(2)},(1,3,2)\bigr)$ and $a_i\mapsto\bigl(\gamma^{(1)}\gamma^{(2)},(1,2,i+1)\bigr)$, $i=3,\dots,n-1$.
\item for the presentation (\ref{codeDn1}), 

$R_1\mapsto\bigl(\gamma^{(1)}\gamma^{(2)},\eA\bigr)$ and $R_i\mapsto\bigl(\gamma^{(1)}\gamma^{(2)},(1,2)(i,i+1)\bigr)$, $i=2,\dots,n-1$.
\item for the presentation (\ref{codeDn2}), 

$r_1\mapsto\bigl(\eCn ,(1,2,3)\bigr)$, $r_2\mapsto\bigl(\gamma^{(1)}\gamma^{(2)},(1,2,3)\bigr)$ and $r_i\mapsto\bigl(\eCn ,(i-1,i,i+1)\bigr)$, $i=3,\dots,n-1$.
\end{itemize}

Let $H$ be the subgroup generated by $a_1,\dots,a_{n-2}$, or by $R_1$, $\dots$, $R_{n-2}$, or by $r_1$, $\dots$, $r_{n-2}$. 

\vskip 2cm
\setlength{\unitlength}{1800sp}
\begingroup\makeatletter\ifx\SetFigFontNFSS\undefined
\gdef\SetFigFontNFSS#1#2#3#4#5{
  \reset@font\fontsize{#1}{#2pt}
  \fontfamily{#3}\fontseries{#4}\fontshape{#5}
  \selectfont}
\fi\endgroup
\begin{picture}(2283,4528)(371,-6560)
{\put(5300,-3100){\circle*{140}}}
\put(100,-2850){\makebox(0,0)[lb]{\smash{{\SetFigFontNFSS{12}{14.4}{\rmdefault}{\mddefault}{\updefault}{\small{$H$}}}}}}
{\put(5300,-3100){\vector(2,3){900}}}
{\put(6300,-1720){\circle*{140}}}
{\put(6300,-4480){\circle*{140}}}
{\put(6250,-4480){\vector(-2,3){900}}}
\put(5150,-3150){\makebox(0,0)[lb]{\smash{{\SetFigFontNFSS{12}{14.4}{\rmdefault}{\mddefault}{\updefault}{\small{$a_3$}}}}}}
{\put(3450,-3100){\circle*{140}}}
{\put(3500,-3050){\vector(2,1){2730}}}
{\put(6230,-4480){\vector(-2,1){2750}}}
\put(3550,-3150){\makebox(0,0)[lb]{\smash{{\SetFigFontNFSS{12}{14.4}{\rmdefault}{\mddefault}{\updefault}{\small{$a_4$}}}}}}
{\put(550,-3100){\circle*{140}}}
{\put(550,-3080){\vector(4,1){5640}}}
{\put(6240,-4500){\vector(-4,1){5630}}}
\put(850,-3150){\makebox(0,0)[lb]{\smash{{\SetFigFontNFSS{12}{14.4}{\rmdefault}{\mddefault}{\updefault}{\small{$a_{n-1}\ \dots\ \dots$}}}}}}
{\put(6210,-1790){\vector(0,-1){2600}}}
\put(5780,-3200){\makebox(0,0)[lb]{\smash{{\SetFigFontNFSS{12}{14.4}{\rmdefault}{\mddefault}{\updefault}{\tiny{$\cdot$}}}}}}
\put(5820,-3200){\makebox(0,0)[lb]{\smash{{\SetFigFontNFSS{12}{14.4}{\rmdefault}{\mddefault}{\updefault}{\tiny{$\cdot$}}}}}}
\put(5860,-3200){\makebox(0,0)[lb]{\smash{{\SetFigFontNFSS{12}{14.4}{\rmdefault}{\mddefault}{\updefault}{\tiny{$\cdot$}}}}}}
{\put(6330,-1790){\vector(0,-1){2600}}}
{\put(6420,-4470){\vector(0,1){2600}}}
{\put(6350,-4480){\vector(3,-2){2020}}}
{\put(8440,-5750){\vector(-1,2){2040}}}
{\put(8430,-5830){\circle*{140}}}
\put(7200,-5150){\makebox(0,0)[lb]{\smash{{\SetFigFontNFSS{12}{14.4}{\rmdefault}{\mddefault}{\updefault}{\small{$a_1$}}}}}}
{\put(6410,-1870){\vector(3,2){1980}}}
{\put(8420,-570){\vector(-1,-2){1950}}}
{\put(8450,-510){\circle*{140}}}
\put(7200,-1300){\makebox(0,0)[lb]{\smash{{\SetFigFontNFSS{12}{14.4}{\rmdefault}{\mddefault}{\updefault}{\small{$a_2$}}}}}}
{\put(8430,-5830){\vector(3,2){1980}}}
{\put(10380,-1840){\vector(-1,-2){1950}}}
{\put(10500,-4450){\circle*{140}}}
\put(8530,-5150){\makebox(0,0)[lb]{\smash{{\SetFigFontNFSS{12}{14.4}{\rmdefault}{\mddefault}{\updefault}{\small{$a_2$}}}}}}
{\put(8490,-460){\vector(3,-2){1980}}}
{\put(10470,-4420){\vector(-1,2){1950}}}
{\put(10500,-1820){\circle*{140}}}
\put(8700,-1300){\makebox(0,0)[lb]{\smash{{\SetFigFontNFSS{12}{14.4}{\rmdefault}{\mddefault}{\updefault}{\small{$a_1$}}}}}}
{\put(10490,-1790){\vector(0,-1){2600}}}
\put(10074,-3200){\makebox(0,0)[lb]{\smash{{\SetFigFontNFSS{12}{14.4}{\rmdefault}{\mddefault}{\updefault}{\tiny{$\cdot$}}}}}}
\put(10105,-3200){\makebox(0,0)[lb]{\smash{{\SetFigFontNFSS{12}{14.4}{\rmdefault}{\mddefault}{\updefault}{\tiny{$\cdot$}}}}}}
\put(10135,-3200){\makebox(0,0)[lb]{\smash{{\SetFigFontNFSS{12}{14.4}{\rmdefault}{\mddefault}{\updefault}{\tiny{$\cdot$}}}}}}
{\put(10610,-1860){\vector(0,-1){2530}}}
{\put(10400,-4470){\vector(0,1){2600}}}
{\put(10620,-4430){\vector(2,3){840}}}
{\put(11520,-3130){\circle*{140}}}
{\put(11450,-3110){\vector(-2,3){840}}}
\put(10500,-3150){\makebox(0,0)[lb]{\smash{{\SetFigFontNFSS{12}{14.4}{\rmdefault}{\mddefault}{\updefault}{\small{$a_3$}}}}}}
{\put(10640,-4460){\vector(2,1){2650}}}
{\put(13300,-3120){\circle*{140}}}
{\put(13230,-3120){\vector(-2,1){2620}}}
\put(12000,-3150){\makebox(0,0)[lb]{\smash{{\SetFigFontNFSS{12}{14.4}{\rmdefault}{\mddefault}{\updefault}{\small{$a_4$}}}}}}
{\put(10700,-4500){\vector(4,1){5450}}}
{\put(16200,-3120){\circle*{140}}}
{\put(16180,-3130){\vector(-4,1){5480}}}
\put(13120,-3150){\makebox(0,0)[lb]{\smash{{\SetFigFontNFSS{12}{14.4}{\rmdefault}{\mddefault}{\updefault}{\small{$\dots\ \ \dots\ \ a_{n-1}$}}}}}}
\put(1400,-6500){\begin{pic}\label{CoxToddaltD3} Coxeter--Todd figure for $(D^+_n,H)$ for the presentation (\ref{codeDn3})
\end{pic}}
\end{picture}

\bigskip
\setlength{\unitlength}{2450sp}
\begingroup\makeatletter\ifx\SetFigFontNFSS\undefined
\gdef\SetFigFontNFSS#1#2#3#4#5{
  \reset@font\fontsize{#1}{#2pt}
  \fontfamily{#3}\fontseries{#4}\fontshape{#5}
  \selectfont}
\fi\endgroup
\begin{picture}(12110,2938)(346,-4855)
{\thinlines
\put(541,-3346){\circle*{144}}}{\put(1486,-3346){\circle*{144}}}{\put(2476,-3346){\circle*{144}}}{\put(10396,-3346){\circle*{144}}}{\put(11386,-3346){\circle*{144}}}{\put(12376,-3346){\circle*{144}}}{\put(3736,-3346){\circle*{144}}}{\put(4726,-3346){\circle*{144}}}{\put(5806,-2266){\circle*{144}}}{\put(5824,-4516){\circle*{144}}}{\put(9181,-3346){\circle*{144}}}{\put(8281,-3346){\circle*{144}}}{\put(7156,-4516){\circle*{144}}}{\put(7111,-2266){\circle*{144}}}{\put(586,-3346){\line(1,0){810}}}{\put(1576,-3346){\line(1,0){810}}}{\put(2566,-3346){\line(1,0){270}}}{\put(2971,-3346){\line(1,0){270}}}{\put(3376,-3346){\line(1,0){270}}}{\put(9226,-3346){\line(1,0){270}}}{\put(9676,-3346){\line(1,0){270}}}{\put(10081,-3346){\line(1,0){270}}}{\put(10486,-3346){\line(1,0){810}}}{\put(11476,-3346){\line(1,0){810}}}{\put(3826,-3346){\line(1,0){810}}}{\put(4748,-3278){\vector(1,1){1035}}}{\put(5761,-4471){\vector(-1,1){1035}}}{\put(8326,-3346){\line(1,0){810}}}{\put(8213,-3414){\vector(-1,-1){1035}}}{\put(5761,-2311){\vector(0,-1){2205}}}{\put(5806,-2311){\line(0,-1){2160}}}{\put(7156,-2311){\line(0,-1){2160}}}{\put(5896,-4426){\line(0,1){2115}}}{\put(7201,-2266){\vector(0,-1){2205}}}{\put(7222,-2332){\vector(1,-1){1035}}}{\put(7066,-2311){\line(0,-1){2160}}}{\put(5851,-2266){\line(1,0){1215}}}{\put(5851,-4516){\line(1,0){1260}}}\put(61,-3706){\makebox(0,0)[lb]{\smash{{\SetFigFontNFSS{12}{14.4}{\rmdefault}{\mddefault}{\updefault}{$H$}}}}}\put(3761,-3256){\makebox(0,0)[lb]{\smash{{\SetFigFontNFSS{12}{14.4}{\rmdefault}{\mddefault}{\updefault}{\small{$R_3$}}}}}}\put(431,-3256){\makebox(0,0)[lb]{\smash{{\SetFigFontNFSS{12}{14.4}{\rmdefault}{\mddefault}{\updefault}{\small{$R_{n-1}$}}}}}}\put(1376,-3256){\makebox(0,0)[lb]{\smash{{\SetFigFontNFSS{12}{14.4}{\rmdefault}{\mddefault}{\updefault}{\small{$R_{n-2}$}}}}}}\put(10331,-3256){\makebox(0,0)[lb]{\smash{{\SetFigFontNFSS{12}{14.4}{\rmdefault}{\mddefault}{\updefault}{\small{$R_{n-2}$}}}}}}\put(11276,-3256){\makebox(0,0)[lb]{\smash{{\SetFigFontNFSS{12}{14.4}{\rmdefault}{\mddefault}{\updefault}{\small{$R_{n-1}$}}}}}}\put(8306,-3256){\makebox(0,0)[lb]{\smash{{\SetFigFontNFSS{12}{14.4}{\rmdefault}{\mddefault}{\updefault}{\small{$R_3$}}}}}}\put(6011,-2176){\makebox(0,0)[lb]{\smash{{\SetFigFontNFSS{12}{14.4}{\rmdefault}{\mddefault}{\updefault}{\small{$R_1$}}}}}}\put(6011,-4786){\makebox(0,0)[lb]{\smash{{\SetFigFontNFSS{12}{14.4}{\rmdefault}{\mddefault}{\updefault}{\small{$R_1$}}}}}}\put(7136,-3391){\makebox(0,0)[lb]{\smash{{\SetFigFontNFSS{12}{14.4}{\rmdefault}{\mddefault}{\updefault}{\small{$R_2$}}}}}}\put(4841,-3391){\makebox(0,0)[lb]{\smash{{\SetFigFontNFSS{12}{14.4}{\rmdefault}{\mddefault}{\updefault}{\small{$R_2$}}}}}}\put(5506,-3436){\makebox(0,0)[lb]{\smash{{\SetFigFontNFSS{12}{14.4}{\rmdefault}{\mddefault}{\updefault}{\small{$R_3,\dots,R_{n-1}$}}}}}}
\put(1700,-5300){\begin{pic}\label{CoxToddaltD1} Coxeter--Todd figure for $(D^+_n,H)$ for the presentation (\ref{codeDn1})
\end{pic}}
\end{picture}

\vskip 1cm
\setlength{\unitlength}{2450sp}
\begingroup\makeatletter\ifx\SetFigFontNFSS\undefined
\gdef\SetFigFontNFSS#1#2#3#4#5{
  \reset@font\fontsize{#1}{#2pt}
  \fontfamily{#3}\fontseries{#4}\fontshape{#5}
  \selectfont}
\fi\endgroup
\begin{picture}(12083,4028)(571,-6260)
{\thinlines
\put(811,-3346){\circle*{144}}}{\put(1711,-5146){\circle*{144}}}{\put(2611,-3346){\circle*{144}}}{\put(3511,-5146){\circle*{144}}}{\put(4024,-3346){\circle*{144}}}{\put(4861,-5146){\circle*{144}}}{\put(4861,-3346){\circle*{144}}}{\put(6391,-2311){\circle*{144}}}{\put(6409,-6181){\circle*{144}}}{\put(7921,-3346){\circle*{144}}}{\put(7939,-5191){\circle*{144}}}{\put(11566,-3346){\circle*{144}}}{\put(10711,-5191){\circle*{144}}}{\put(12574,-5191){\circle*{144}}}{\put(8866,-5191){\circle*{144}}}{\put(9721,-3346){\circle*{144}}}{\put(1675,-5119){\vector(-1,2){855}}}{\put(2548,-3382){\vector(-1,-2){855}}}{\put(2620,-3391){\vector(-1,-2){855}}}{\put(3520,-5119){\vector(-1,2){890}}}{\put(2746,-3346){\line(1,0){270}}}{\put(3196,-3346){\line(1,0){270}}}{\put(3646,-3346){\line(1,0){270}}}{\put(3646,-5146){\line(1,0){270}}}{\put(4051,-5146){\line(1,0){270}}}{\put(4501,-5146){\line(1,0){270}}}{\put(4879,-5137){\vector(-1,2){890}}}{\put(4861,-3391){\vector(0,-1){1710}}}{\put(4096,-3346){\vector(1,0){720}}}{\put(6391,-2311){\vector(-3,-2){1485}}}{\put(4951,-5146){\vector(1,2){1395}}}{\put(4861,-5146){\vector(3,-2){1485}}}{\put(6346,-6136){\vector(-1,2){1395}}}{\put(7921,-3301){\vector(-3,2){1485}}}{\put(7921,-5146){\vector(0,1){1710}}}{\put(6481,-2401){\vector(1,-2){1395}}}{\put(7831,-3346){\vector(-1,-2){1395}}}{\put(6467,-6160){\vector(3,2){1485}}}{\put(10693,-5092){\vector(1,2){855}}}{\put(11647,-3418){\vector(1,-2){855}}}{\put(9856,-5191){\line(1,0){270}}}{\put(10306,-5191){\line(1,0){270}}}{\put(11116,-3346){\line(1,0){270}}}{\put(10711,-3346){\line(1,0){270}}}{\put(8757,-5191){\vector(-1,0){720}}}{\put(8794,-5182){\vector(-1,2){890}}}{\put(9856,-3346){\line(1,0){270}}}{\put(10306,-3346){\line(1,0){270}}}{\put(9451,-5191){\line(1,0){270}}}{\put(9046,-5191){\line(1,0){270}}}{\put(9721,-3436){\vector(-1,-2){855}}}{\put(8840,-5146){\vector(-1,2){890}}}{\put(901,-3346){\vector(1,0){1665}}}{\put(1801,-5146){\vector(1,0){1665}}}{\put(12466,-5191){\vector(-1,0){1665}}}{\put(8011,-3346){\vector(1,0){1665}}}\put(386,-3166){\makebox(0,0)[lb]{\smash{{\SetFigFontNFSS{12}{14.4}{\rmdefault}{\mddefault}{\updefault}{$H$}}}}}\put(1106,-3976){\makebox(0,0)[lb]{\smash{{\SetFigFontNFSS{12}{14.4}{\rmdefault}{\mddefault}{\updefault}{$r_{n-1}$}}}}}\put(2096,-4651){\makebox(0,0)[lb]{\smash{{\SetFigFontNFSS{12}{14.4}{\rmdefault}{\mddefault}{\updefault}{$r_{n-2}$}}}}}\put(4166,-3976){\makebox(0,0)[lb]{\smash{{\SetFigFontNFSS{12}{14.4}{\rmdefault}{\mddefault}{\updefault}{$r_3$}}}}}\put(5156,-3256){\makebox(0,0)[lb]{\smash{{\SetFigFontNFSS{12}{14.4}{\rmdefault}{\mddefault}{\updefault}{$r_2$}}}}}\put(5156,-5416){\makebox(0,0)[lb]{\smash{{\SetFigFontNFSS{12}{14.4}{\rmdefault}{\mddefault}{\updefault}{$r_1$}}}}}\put(6731,-5461){\makebox(0,0)[lb]{\smash{{\SetFigFontNFSS{12}{14.4}{\rmdefault}{\mddefault}{\updefault}{$r_2$}}}}}\put(6821,-3256){\makebox(0,0)[lb]{\smash{{\SetFigFontNFSS{12}{14.4}{\rmdefault}{\mddefault}{\updefault}{$r_1$}}}}}\put(11006,-4741){\makebox(0,0)[lb]{\smash{{\SetFigFontNFSS{12}{14.4}{\rmdefault}{\mddefault}{\updefault}{$r_{n-1}$}}}}}\put(7811,-4786){\makebox(0,0)[lb]{\smash{{\SetFigFontNFSS{12}{14.4}{\rmdefault}{\mddefault}{\updefault}{$r_3$}}}}}\put(8441,-3931){\makebox(0,0)[lb]{\smash{{\SetFigFontNFSS{12}{14.4}{\rmdefault}{\mddefault}{\updefault}{$r_4$}}}}}
\put(1500,-6700){\begin{pic}\label{CoxToddaltD2a} Coxeter--Todd figure for $(D^+_n,H)$ for the presentation (\ref{codeDn2}), $n$ even
\end{pic}}
\end{picture}

\vskip 1cm
\setlength{\unitlength}{2450sp}
\begingroup\makeatletter\ifx\SetFigFontNFSS\undefined
\gdef\SetFigFontNFSS#1#2#3#4#5{
  \reset@font\fontsize{#1}{#2pt}
  \fontfamily{#3}\fontseries{#4}\fontshape{#5}
  \selectfont}
\fi\endgroup
\begin{picture}(12283,4528)(371,-6560)
{\thinlines
\put(4861,-5146){\circle*{144}}}{\put(4861,-3346){\circle*{144}}}{\put(6391,-2311){\circle*{144}}}{\put(6409,-6181){\circle*{144}}}{\put(7921,-3346){\circle*{144}}}{\put(7939,-5191){\circle*{144}}}{\put(11566,-3346){\circle*{144}}}{\put(10711,-5191){\circle*{144}}}{\put(12574,-5191){\circle*{144}}}{\put(4006,-5146){\circle*{144}}}{\put(451,-3346){\circle*{144}}}{\put(2251,-3346){\circle*{144}}}{\put(1351,-5146){\circle*{144}}}{\put(3151,-5146){\circle*{144}}}{\put(8821,-3346){\circle*{144}}}{\put(9739,-5191){\circle*{144}}}{\put(4861,-5101){\vector(0,1){1710}}}{\put(4906,-3301){\vector(3,2){1485}}}{\put(6346,-2356){\vector(-1,-2){1395}}}{\put(6346,-6136){\vector(-3,2){1485}}}{\put(4951,-3346){\vector(1,-2){1395}}}{\put(6436,-2311){\vector(3,-2){1485}}}{\put(7921,-3436){\vector(0,-1){1710}}}{\put(7876,-5191){\vector(-1,2){1410}}}{\put(6436,-6136){\vector(1,2){1395}}}{\put(7952,-5170){\vector(-3,-2){1485}}}{\put(10693,-5092){\vector(1,2){855}}}{\put(11647,-3418){\vector(1,-2){855}}}{\put(11116,-3346){\line(1,0){270}}}{\put(10261,-3346){\line(1,0){270}}}{\put(10711,-3346){\line(1,0){270}}}{\put(4852,-3409){\vector(-1,-2){855}}}{\put(4096,-5146){\vector(1,0){720}}}{\put(1342,-5128){\vector(-1,2){855}}}{\put(2260,-3409){\vector(-1,-2){855}}}{\put(2206,-3346){\vector(-1,-2){855}}}{\put(3160,-5119){\vector(-1,2){855}}}{\put(2341,-3346){\line(1,0){270}}}{\put(2791,-3346){\line(1,0){270}}}{\put(3241,-3346){\line(1,0){270}}}{\put(3241,-5146){\line(1,0){270}}}{\put(3601,-5146){\line(1,0){270}}}{\put(3646,-3346){\line(1,0){270}}}{\put(4096,-3346){\line(1,0){270}}}{\put(4501,-3346){\line(1,0){270}}}{\put(8733,-3346){\vector(-1,0){719}}}{\put(7876,-5146){\vector(1,2){890}}}{\put(7939,-5200){\vector(1,2){890}}}{\put(9685,-5119){\vector(-1,2){855}}}{\put(9901,-5191){\line(1,0){270}}}{\put(10306,-5191){\line(1,0){270}}}{\put(9001,-3346){\line(1,0){270}}}{\put(9406,-3346){\line(1,0){270}}}{\put(9856,-3346){\line(1,0){270}}}{\put(12466,-5191){\vector(-1,0){1665}}}{\put(541,-3346){\vector(1,0){1665}}}{\put(1441,-5146){\vector(1,0){1665}}}{\put(9676,-5191){\vector(-1,0){1665}}}\put(86,-3166){\makebox(0,0)[lb]{\smash{{\SetFigFontNFSS{12}{14.4}{\rmdefault}{\mddefault}{\updefault}{$H$}}}}}\put(5256,-3256){\makebox(0,0)[lb]{\smash{{\SetFigFontNFSS{12}{14.4}{\rmdefault}{\mddefault}{\updefault}{$r_2$}}}}}\put(5256,-5416){\makebox(0,0)[lb]{\smash{{\SetFigFontNFSS{12}{14.4}{\rmdefault}{\mddefault}{\updefault}{$r_1$}}}}}\put(6831,-5461){\makebox(0,0)[lb]{\smash{{\SetFigFontNFSS{12}{14.4}{\rmdefault}{\mddefault}{\updefault}{$r_2$}}}}}\put(6921,-3256){\makebox(0,0)[lb]{\smash{{\SetFigFontNFSS{12}{14.4}{\rmdefault}{\mddefault}{\updefault}{$r_1$}}}}}\put(11106,-4741){\makebox(0,0)[lb]{\smash{{\SetFigFontNFSS{12}{14.4}{\rmdefault}{\mddefault}{\updefault}{$r_{n-1}$}}}}}\put(756,-3976){\makebox(0,0)[lb]{\smash{{\SetFigFontNFSS{12}{14.4}{\rmdefault}{\mddefault}{\updefault}{$r_{n-1}$}}}}}\put(1780,-4651){\makebox(0,0)[lb]{\smash{{\SetFigFontNFSS{12}{14.4}{\rmdefault}{\mddefault}{\updefault}{$r_{n-2}$}}}}}\put(4221,-4696){\makebox(0,0)[lb]{\smash{{\SetFigFontNFSS{12}{14.4}{\rmdefault}{\mddefault}{\updefault}{$r_3$}}}}}\put(7866,-3931){\makebox(0,0)[lb]{\smash{{\SetFigFontNFSS{12}{14.4}{\rmdefault}{\mddefault}{\updefault}{$r_3$}}}}}\put(8451,-4786){\makebox(0,0)[lb]{\smash{{\SetFigFontNFSS{12}{14.4}{\rmdefault}{\mddefault}{\updefault}{$r_4$}}}}}
\put(1400,-6650){\begin{pic}\label{CoxToddaltD2b} Coxeter--Todd figure for $(D^+_n,H)$ for the presentation (\ref{codeDn2}), $n$ odd
\end{pic}}
\end{picture}

\vskip .4cm
The  Coxeter--Todd figures for the presentations  (\ref{codeDn3}), (\ref{codeDn1}) and  (\ref{codeDn2}) provide three normal forms for elements of the alternating group $D^+_n$. We omit details.

\paragraph{Acknowledgements.} We thank I. Marin, J. Michel, O. Brunat, C. Blanchet for useful discussions.

\end{document}